\documentclass[
	,12pt%
	,pagesize%
	,headings=small%
	,paper=a4%
	,parskip=false%
	,abstract=true
	,toc=bibliography
	,DIV=12%
]{scrartcl}

\newcommand{\showbrief}{}

\addtokomafont{sectioning}{\normalfont\bfseries}
\setkomafont{title}{\normalfont}
\setkomafont{subtitle}{\normalfont}
\setkomafont{author}{\normalfont}
\setkomafont{date}{\normalfont}
\addtokomafont{pageheadfoot}{\scshape\small}
\setkomafont{caption}{\footnotesize}
\setkomafont{captionlabel}{\usekomafont{caption}\bfseries}
\setcapindent{0pt}

\usepackage{xspace}

\usepackage[utf8]{inputenc}
\usepackage[T1]{fontenc}
\usepackage{csquotes}
\usepackage[english]{babel}
\usepackage{graphicx}
\graphicspath{{Pictures/}{/}}
\usepackage[export]{adjustbox}	

\usepackage[dvipsnames]{xcolor}
\usepackage{tikz-cd}

\usepackage[final]{fixme}
\usepackage[nointlimits]{amsmath}
\usepackage{amssymb,mathrsfs,mathtools}
\usepackage{newtxtext}
\usepackage{newtxmath}
\AtBeginDocument{%
  \mathchardef\standardeq=\mathcode`=
  \mathchardef\standardless=\mathcode`<
  \mathchardef\standardgreater=\mathcode`>
  \mathcode`="8000
  \mathcode`<"8000
  \mathcode`>"8000
}
\begingroup\lccode`~`=\lowercase{\endgroup
  \def~}{\mathrel{\mspace{0.936mu}\standardeq\mspace{0.936mu}}}
\begingroup\lccode`~`<\lowercase{\endgroup
  \def~}{\mathrel{\mspace{0.153mu}\standardless\mspace{0.153mu}}}
\begingroup\lccode`~`>\lowercase{\endgroup
  \def~}{\mathrel{\mspace{0.153mu}\standardgreater\mspace{0.153mu}}}

\usepackage{upgreek}
\usepackage{braket}
\usepackage{cancel}
\usepackage{siunitx}
\usepackage{aliascnt} 
\usepackage[amsmath,thmmarks]{ntheorem}
\usepackage{xkeyval}
\usepackage{authblk}
\usepackage[final,%
	bookmarks,%
	bookmarksdepth=3,%
	breaklinks=true,%
	colorlinks=true,%
	urlcolor=NavyBlue,%
	linkcolor=NavyBlue,%
	citecolor=ForestGreen,%
]{hyperref}%
\usepackage[all]{hypcap}

\AtBeginEnvironment{abstract}{\footnotesize}
\makeatletter
\patchcmd{\@maketitle}{\huge}{\Large}{}{}
\makeatother

\newcommand{\Curve}{\gamma}

\DeclareDocumentCommand{\Path}{ o }{
	\IfValueTF{#1}{%
		\varGamma_{\!#1}%
	}{%
		\varGamma%
	}%
}

\DeclareDocumentCommand{\dotPath}{ o }{
	\IfValueTF{#1}{%
		\dot{\varGamma}_{\!#1}%
	}{%
		{\dot{\varGamma}}%
	}%
}

\newcommand{\Mfld}{H^s_{\mathrm{i,r}}(\Domain, \AmbSpace)}

\DeclareMathOperator{\BiLip}{BiLip}

\DeclareDocumentCommand{\dist}{ o o o }{
	\IfValueTF{#1}{
		\text{dist}_{#1}\IfValueTF{#2}{(#2,#3)}{}
	}{
		\text{dist}\IfValueTF{#2}{(#2,#3)}{}
	}
}

\DeclareDocumentCommand{\distC}{ o o }{
	\varrho_{\gamma}\IfValueTF{#1}{(#1,#2)}{}
}

\DeclareDocumentCommand{\LineEl}{ o o }{
	\IfValueTF{#1}{
	  \omega_{#1}\IfValueTF{#2}{(#2)}{}
	}{
	  \omega\IfValueTF{#2}{(#2)}{}
	}
}

\DeclareDocumentCommand{\dLineEl}{ o o }{
	\IfValueTF{#1}{
	  \dd \omega_{#1}\IfValueTF{#2}{(#2)}{}
	}{
    \dd \omega\IfValueTF{#2}{(#2)}{}
	}
}

\DeclareDocumentCommand{\LineElC}{ o }{
	\IfValueTF{#1}{
        \omega_\Curve(#1)
	}{
        \omega_\Curve
	}
}

\DeclareDocumentCommand{\dLineElC}{ o }{
	\IfValueTF{#1}{
	  \dd \omega_\Curve(#1)
	}{
        \dd \omega_\Curve
	}
}

\DeclareDocumentCommand{\LebesgueM}{ o }{
	\IfValueTF{#1}{
        \lambda(#1)
	}{
        \lambda
	}
}

\DeclareDocumentCommand{\dLebesgueM}{ o }{
	\IfValueTF{#1}{
	  \dd \lambda(#1)
	}{
        \dd \lambda
	}
}

\DeclareDocumentCommand{\Speed}{ o }{
	\IfValueTF{#1}{
	  	h_{#1}
	}{
      	h
	}
}

\DeclareDocumentCommand{\InvSpeed}{ o }{
	\IfValueTF{#1}{
	  	H_{#1}
	}{
      	H
	}
}

\newcommand{\brief}[1]{\ifthenelse{\isundefined{\showbrief}}{}{{\color{NavyBlue}{\bigskip\emph{Brief:} #1\newline}}}}%

\newcommand{\Op}[2]{\mathcal{R}_{#1}^{#2}}

\newcommand{\Circle}{\mathbb{T}}

\newcommand{\TP}{\mathrm{TP}}

\newcommand{\AmbDim}{n} 
\newcommand{\AmbSpace}{{\R^\AmbDim}}
\newcommand{\Domain}{{\Circle}}

\DeclareDocumentCommand{\Hess}{ O{} }{\operatorname{Hess}_{#1}}

\DeclareDocumentCommand{\converges}{ o }{
	\mathbin{%
		\IfValueTF{#1}{%
			\mathrel{\vbox{\offinterlineskip\ialign{%
				\hfil##\hfil\cr
				$\scriptscriptstyle#1$\cr
				$-\!\!\!-\!\!\!\rightarrow$\cr
			}}}
		}{%
			-\!\!\!-\!\!\!\rightarrow
		}%
	}%
}

\DeclareDocumentCommand{\wconverges}{ o }{
	\mathbin{%
		\IfValueTF{#1}{%
			\mathrel{\vbox{\offinterlineskip\ialign{%
				\hfil##\hfil\cr
				$\scriptscriptstyle#1$\cr
				$-\!\!\!-\!\!\!\rightharpoonup$\cr
			}}}
		}{%
			-\!\!\!-\!\!\!\rightharpoonup
		}%
	}%
}

\newcommand{\dd}{\mathop{}\!\mathrm{d}}

\newcommand{\ceq}{\coloneqq}

\newcommand{\R}{{\mathbb{R}}}

\newcommand{\N}{\mathbb{N}}
\newcommand{\Z}{{\mathbb{Z}}}

\DeclarePairedDelimiterXPP{\pars}[1]{\mathop{}}{\lparen}{\rparen}{}{#1}
\DeclarePairedDelimiterXPP{\abs}[1]{\mathop{}}{\lvert}{\rvert}{}{#1}
\DeclarePairedDelimiterXPP{\norm}[1]{\mathop{}}{\lVert}{\rVert}{}{#1}
\DeclarePairedDelimiterXPP{\seminorm}[1]{\mathop{}}{\lbrack}{\rbrack}{}{#1}
\DeclarePairedDelimiterXPP{\inner}[1]{\mathop{}}{\langle}{\rangle}{}{#1}
\DeclarePairedDelimiterXPP{\iinner}[1]{\mathop{}}{\langle\!\langle}{\rangle\!\rangle}{}{#1}
\DeclarePairedDelimiterXPP{\brackets}[1]{\mathop{}}{\lbrack}{\rbrack}{}{#1}
\DeclarePairedDelimiterXPP{\braces}[1]{\mathop{}}{\lbrace}{\rbrace}{}{#1}

\DeclarePairedDelimiterXPP{\floor}[1]{\mathop{}}{\lfloor}{\rfloor}{}{#1}
\DeclarePairedDelimiterXPP{\ceil}[1]{\mathop{}}{\lceil}{\rceil}{}{#1}

\DeclarePairedDelimiterXPP{\intervalcc}[1]{\mathop{}}{\lbrack}{\rbrack}{}{#1}
\DeclarePairedDelimiterXPP{\intervalco}[1]{\mathop{}}{\lbrack}{\rparen}{}{#1}
\DeclarePairedDelimiterXPP{\intervaloc}[1]{\mathop{}}{\lparen}{\rbrack}{}{#1}
\DeclarePairedDelimiterXPP{\intervaloo}[1]{\mathop{}}{\lparen}{\rparen}{}{#1}

\DeclarePairedDelimiterXPP{\myset}[2]{\mathop{}}{\lbrace}{\rbrace}{}{#1\,\delimsize\vert\,\mathopen{}#2}

\newcommand{\nocontentsline}[3]{}
\newcommand{\tocless}[2]{\bgroup\let\addcontentsline=\nocontentsline#1{#2}\egroup}

\newcommand{\AL}{\mathcal{A}}

\newcommand{\D}{\mathrm{D}}
\newcommand{\di}{\mathrm{d}}
\newcommand{\TPpq}[1][p,q]{\ensuremath{\TP^{\left(#1\right)}}}
\newcommand{\TPp}{\TPpq[p,2]}

\newcommand{\T}{\mathrm{T}}

\newcommand{\Bilip}{\mathrm{BiLip}}
\newcommand{\pdist}{\mathrm{dist}}
\newcommand{\EL}{\mathcal{E}}
\newcommand{\HL}{\mathcal{H}}
\newcommand{\UL}{\mathcal{U}}
\newcommand{\OL}{\mathcal{O}}
\newcommand{\ML}{\mathcal{M}}

\newenvironment{beweis}[1][\unskip]{
		\textbf{Proof #1:}
	}{
		\hbox{}
		\newline
		\hbox{}
		\hfill\includegraphics[scale=0.05]{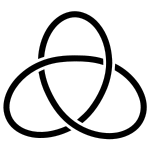}
	}

\DeclareDocumentCommand{\Graph}{ O{} O{} o o}{
	\IfValueTF{#3}{
	  \IfValueTF{#4}{
	  	\mathrm{Graph}^{#1}_{#2}(#3;#4)
	  }{
		\mathrm{Graph}^{#1}_{#2}(#3)
	  }
	}{
		\mathrm{Graph}^{#1}_{#2}
	}
}

\DeclareDocumentCommand{\Emb}{ O{} O{} o o}{
	\IfValueTF{#3}{
	  \IfValueTF{#4}{
	  	\mathrm{Emb}^{#1}_{#2}(#3;#4)
	  }{
		\mathrm{Emb}^{#1}_{#2}(#3)
	  }
	}{
		\mathrm{Emb}^{#1}_{#2}
	}
}

\DeclareDocumentCommand{\Sobo}{ O{} O{} o o}{
	\IfValueTF{#3}{
	  \IfValueTF{#4}{
	  	W^{#1}_{#2}(#3;#4)
	  }{
	  	W^{#1}_{#2}(#3)
	  }
	}{
	  W^{#1}_{#2}
	}
}

\DeclareDocumentCommand{\Bessel}{ O{} O{} o o}{
	\IfValueTF{#3}{
	  \IfValueTF{#4}{
	  	H^{#1}_{#2}(#3;#4)
	  }{
	  	H^{#1}_{#2}(#3)
	  }
	}{
	  H^{#1}_{#2}
	}
}

\DeclareDocumentCommand{\Holder}{ O{} O{} o o}{
	\IfValueTF{#3}{
	  \IfValueTF{#4}{
	  	C^{#1}_{#2}(#3;#4)
	  }{
	  	C^{#1}_{#2}(#3)
	  }
	}{
	  C^{#1}_{#2}
	}
}
\DeclareDocumentCommand{\HolderC}{ O{} O{} }{\Holder[#1][#2][\Circle][\AmbSpace]}

\DeclareDocumentCommand{\Lebesgue}{ O{} O{} o o}{
	\IfValueTF{#3}{
	  \IfValueTF{#4}{
	  	L^{#1}_{#2}(#3;#4)
	  }{
	  	L^{#1}_{#2}(#3)
	  }
	}{
	  L^{#1}_{#2}
	}
}
\DeclareDocumentCommand{\LebesgueC}{ O{} O{} }{\Lebesgue[#1][#2][\Circle][\AmbSpace]}

\usepackage[capitalise,nameinlink]{cleveref}

\newcommand{\aref}[1]{\cref{#1}}

\newtheorem{theorem}{Theorem}[section]
\newtheorem{lemma}[theorem]{Lemma}

\newtheorem{corollary}[theorem]{Corollary}

\theoremstyle{break}

\theoremstyle{plain}
\theorembodyfont{\normalfont}
\newtheorem{definition}[theorem]{Definition}
\newtheorem{example}[theorem]{Example}
{%
    \theoremsymbol{\ensuremath{\Diamond}}%
    \newtheorem{remark}[theorem]{Remark}%
}

\theoremstyle{break}

\theoremheaderfont{\itshape}
\theorembodyfont{\upshape}
\theoremstyle{nonumberplain}
\theoremseparator{.}
\theoremsymbol{\ensuremath{\Box}}

\usepackage[backend=biber,backref, giveninits=true,isbn=false,maxbibnames=50,style=alphabetic]{biblatex}
\title{Convergence of gradient flows on knotted curves}

\author[1]{Elias Döhrer}
\author[2]{Nicolas Freches}

\affil[1]{Chemnitz University of Technology, Chemnitz, Germany}
\affil[2]{RWTH Aachen University}

\begin{document}

\maketitle
We prove full convergence of gradient-flows of the arc-length restricted tangent point energies in the Hilbert-case towards critical points.
This is done through a \textit{\L{}ojasiewicz-Simon gradient inequality} for these energies. In order to do so, we prove, that the tangent-point energies are anlytic on the manifold of immersed embeddings and that their Hessian is Fredholm with index zero on the manifold of arc-length parametrized curves.
As a by-product, we also show that the metric on this manifold defined by the first author in \cite{DohrerReiterSchumacherCompleteRiemannianMetric2025} is analytic. 
 \tableofcontents
\section{Introduction}
In order to find appealing representatives within a knot class, 
self-repulsive energies have been investigated over the last forty years.
These energies also prove to be quite useful when modeling and simulating topological effects in physical processes.
This line of research started when Fukuhara introduced the energy of a polygonal knot, motivated by a Coulomb potential (see \cite{FukuharaENERGYKNOT1988}).
O'Hara extended this approach and defined  a family of repulsive energies in \cite{oharaEnergyKnot1991}, \cite{oharaFamilyOfKnotenergies92} and \cite{oharaEnergyOfKnots2-94}.
For a closed, regular, Lipschitz curve $\gamma:\Domain:= \R/\Z \cong \mathbb{S}^1 \rightarrow \AmbSpace$ and $\alpha>0,p>1$, these energies are given by   
\begin{equation}\label{equ:DefOHaraE}
    \nonumber
    \mathscr{O}^{\alpha,p}(\gamma)
    \ceq
        \int_{\R/\Z}
        \int_{\R/\Z}
            \left(
                \frac{1}{\abs{\gamma(y)-\gamma(x)}^\alpha}
                -
                \frac{1}{\dist[\gamma](y,x)^\alpha}
            \right)^{p/2}
        \abs{\gamma'(x)}
        \abs{\gamma'(y)}
        \dd y  
        \dd x
        ,
\end{equation}
where $\dist[\gamma](y,x)$ denotes the distance of $\gamma(y)$ and $\gamma(x)$ on $\gamma(\R/\Z)$.
In this article, we consider the \textit{tangent-point energies}, which first appear in a symmetrized version for the special case $ q=2 $ in the work of Buck and Orloff \cite{BuckOrloff95}.
Gonzalez and Maddocks suggested in \cite{GonzalezMaddocksGlobalCurvatureThickness99} a whole family of tangent-point energies.
The tangent-point energy of a $C^1$-immersion $\gamma: \Domain\rightarrow\AmbSpace$ with parameter $q\in[1,\infty)$ is given by
\begin{equation}\label{equ:DefTPEGeometric}
    \nonumber
    \mathrm{TP}^q(\gamma)
        \ceq
        \int_\Domain
            \int_\Domain
                \frac{1}{r_{\mathrm{TP}(\gamma)}(x,y)^q}
            \abs{\gamma'(y)}
            \abs{\gamma'(x)}
        \dd y   
        \dd x.
\end{equation}
In the above line, 
\[
\frac{1}{r_{\mathrm{TP}(\gamma)}(x,y)}
    \ceq
        2
            \frac{
                \abs{
                    \gamma(y)-\gamma(x)
                    - 
                    \frac{\gamma'(x)}{\abs{\gamma'(x)}}
                    \inner{\frac{\gamma'(x)}{\abs{\gamma'(x)}}, \gamma(y)-\gamma(x)}
                }
            }
            {\abs{\gamma(y)-\gamma(x)}^2}
\]
denotes the reciprocal radius of the smallest circle passing through $\gamma(x)$ and $\gamma(y)$ while also being tangent to $\gamma'(x)$ at $\gamma(x)$.

A self-avoidance property holds true for every $q>2$.
One can show, that every closed, rectifiable curve with finite length and finite tangent-point energy has to be embedded, see \cite[Theorem 1.1]{strzeleckiTangentpointSelfavoidanceEnergies2012}. In the same paper, Strzelecki and von der Mosel showed a regularizing property of the tangent point energy, namely, that any arc-length parametrized injective curve $ \gamma:\Domain\to\R^n  $, with finite energy, is of class $ C^{1,1-\frac2q} $, see \cite[Theorem 1.3]{strzeleckiTangentpointSelfavoidanceEnergies2012}.
The tangent-point energies can also be generalized for "admissible" $k$-dimensional subsets of $ \R^n $,  see \cite{strzeleckiTangentPointRepulsivePotentials2013}.
This generalization still exhibits self-avoidance and regularizing properties for $q>2k$.

In \cite{blattEnergySpacesTangent2013}, Blatt proved that the natural energy space of the tangent-point energy is the fractional Sobolev Slobodeckij space $W^{2-k/q, q}$, that means, that any arclength parametrized injective $ C^1 $-curve $ \gamma $ has finite energy, iff it belongs to this space. 
The assumption q>2 disallows the use of Hilbert space theory. 
Therefore, Blatt and Reiter introduced to \textit{generalized tangent-point energies} in \cite{blattRegularityTheoryTangentpoint2015}.
For a $ C^1 $-immersion $ \gamma:\Domain\to\R^n$ and $ p,q\in[1,\infty) $ we define
\begin{equation}\label{equ:DefTPEGeneralized}
    \begin{split}
    \mathrm{TP}^{(p,q)}
        &\ceq
        2^{-q}
        \int_\Domain
        \int_\Domain
            \frac{1}{r_{\mathrm{TP}(\gamma)}(x,y)^q}
            \frac{1}{\abs{\gamma(y)-\gamma(x)}^{p-2q}}
        \abs{\gamma'(x)}
        \abs{\gamma'(y)}
        \dd y   
        \dd x
        \\
        &=
            \int_\Domain
            \int_\Domain
                \frac{
                    \abs{P^\perp_{\gamma'(x)}(\gamma(y)-\gamma(x))}^q
                }
                {
                    \abs{\gamma(y)-\gamma(x)}^p
                }
            \abs{\gamma'(x)}
            \abs{\gamma'(y)}
            \dd y   
            \dd x
        \quad ,
    \end{split}
\end{equation}
where $P^\perp_{\gamma'(x)}$ denotes the orthogonal projector onto the orthogonal complement of $\R \gamma'(x)$.
For $q>1$ and $p \in \intervaloo{q+2,2q+1}$, these energies still exhibit the self-avoidance property. 
More percisely, one can bound the \textit{Gromov-distorsion} of absolutely continuous, embedded curves $\gamma$ solely in terms of the energy, see \cite[Proposition 2.7]{blatt:2014:RegularityTheoryTangentpoint}.
This means, that there is a $C(p,q) \in (0,\infty)$ such that for an absolutely continuous, embedded curve $\gamma$ 
\begin{equation}\label{equ:DefDistor}
    \mathrm{distor}(\gamma)
    \ceq 
    \sup_{x\neq y}
        \frac{\dist[\gamma](y,x)}{\abs{\gamma(y)-\gamma(x)}}
    \leq
    C(p,q) 
    \mathrm{TP}^{q,p}(\gamma).
\end{equation}
Here $\dist[\gamma](x,y)$ denotes the intrinsic distance of $\gamma(x),\gamma(y) \in \gamma(\Domain)$.
If $\gamma$ is parametrized arc-length, this coincides with the Bi-Lipschitz constant.
This constant is the smallest constant $C\ge 1$ such that
\[
    \frac{1}{C}\dist[\Domain](x,y) 
    \leq \abs{\gamma(x)-\gamma(y)}
    \leq
    C \dist[\Domain](x,y) 
    ,
\]
where $\dist[\Domain](x,y)$ denotes the distance $x$ and $y$ on $\Domain$.
The scale-invariant case ($p=q+2$) poses some difficulties and requires more sophisticated techniques, see \cite{BlattReiterSchikorraVorderobermeierScaleinvariantTangentpointEnergies2024}.

Furthermore, Blatt and Reiter show in \cite[Theorem 1.1]{blattRegularityTheoryTangentpoint2015}, that curves, parametrized by constant speed belong to the space $W^{(p-1)/q,q}$, which, for $ q=2 $, is the Hilbert space $H^s:=W^{s,2},\; s\in \intervaloo{\frac{3}{2},2}$, $ s=\frac{p-1}{2} $. 
In the same paper, smoothness of critical points under a fixed length constraint has been established, see \cite[Theorem 1.5]{blattRegularityTheoryTangentpoint2015}.

In \cite{DohrerReiterSchumacherCompleteRiemannianMetric2025}, Reiter and Schumacher together with the first author exploited the self-avoidance and the characterization of the energy spaces to design a Riemannian metric on the manifold of injective, regular $W^{s,2}$-curves. 
This metric is inspired by the generalized tangent-point energies.
A by-product of the analysis of the Riemannian structure is the smoothness of the functionals $\mathrm{TP}^{(2s+1,2)}$ for $s\in \intervaloo{\frac{3}{2},2}$.

The definition of a knot energy was, at least partially, 
motivated by the idea to disentangle embedded curves by means of gradient flows that respect their topology.
Over the years, a couple of results in this direction appeared.
In \cite{BlattGradientFlowHara2018}, Blatt considered the $L^2$-gradient flow of Ohara's knot energies and was able to show short- and longtime existence and strong convergence after reparametrization using a \L{}ojasiewicz--Simon gradient inequality. See \cite{HeEulerLagrangeEquAndHeatFlowOfMoebiusEnergy2000}, 
\cite{BlattGradFlowMoebiusNearMinimizers} or
\cite{Blatt2020GradFlowMoebius-epsilonRegularity} for similar approaches to the more challenging problem of the $ L^2 $-flow of the Möbius energy.
In \cite{ReiterSchumacherSobolevGradientsObius2021} a Sobolev gradient flow was considered for the Möbius energy $\mathscr{O}^{2,2}$.
There, short time existence has been established in the space $H^{\frac{3}{2}+\varepsilon}$.
A similar approach was considered in \cite{knappmannSpeedPreservingHilbert2022} for a gradient flow of the integral Menger curvature projected onto the nullspace of the logarithmic strain and thereby controlling the parametrization along the flow. 
Short- and longtime existence were established as well as weak subconvergence. 
In \cite{MattSteenebruggeVonDerMoselBanachGradientFlows2023} a minimizing movement approach for the generalized tangent point energies, O'hara energies and integral Menger curvatures in the Banach case has been investigated.

Lastly in \cite{FrechesSchumacherSteenebruggevonderMoselPalaisSmaleConditionGeometric2025} the gradient flow for $\TPp$ was investigated on the submanifold of curves parametrized by arclength. 
There, short- and longtime existence are established as well as strong subconvergence to a critical point by means of the Palais--Smale condition. 
In that regard, \cite{FrechesSchumacherSteenebruggevonderMoselPalaisSmaleConditionGeometric2025} follows a similar approach as \cite{OkabeSchraderConvergenceSobolevGradient2023}, where comparable methods were used to establish results for the elastic bending energy.

We further follow Okabe and Schrader's approach in order to strengthen the results of \cite{FrechesSchumacherSteenebruggevonderMoselPalaisSmaleConditionGeometric2025} to obtain the full convergence of the aforementioned flow, by means of a \L{}ojasiewicz--Simon gradient inequality. 

In order to establish a \L{}ojasiewicz--Simon gradient inequality, it suffices to prove two facts.
First, we show that the tangent-point energy is real analytic on the manifold of immersed $ H^s $-embeddings. We secondly prove that the Hessian, at a critical point, induces a Fredholm Operator of index zero. 
These results yields the \L{}ojasiewicz-Simon inequality and, in combination with the analytic structure of the manifold of $W^{s,2}$-embeddings parametrized by arc length, suffice to show the desired convergence result by standard arguments. 
For the precise statements of our results we defined the manifolds:
\begin{align*}
     \Mfld&\ceq
    \{
        \gamma \in \Bessel[s](\Domain, \AmbSpace)
        \vert
        \gamma \text{ is injective and }
        \gamma \text{ is regular}
    \},
    \\
    \AL^s&\ceq \left\{\gamma\in \Mfld\vert \left|\gamma'(x) \right|=1\text{ for all } x\in \Domain,\gamma \text{ is injective}\right\},\\
    \AL^s_0&\ceq \AL^s \cap \{ \gamma(0)=0\}.
\end{align*}
\begin{theorem}\label{thm:analyticityTPp}
    Let $s=\frac{p-1}{2}\in \intervaloo{\frac{3}{2},2}$.
    The tangent-point energy $\TPp:\Mfld \rightarrow \R$, defined in \eqref{equ:DefTPEGeneralized},  is real analytic.
\end{theorem}
As a consequence, we also obtain the following Corollary.
\begin{corollary}
    Let $s\in \intervaloo{\frac{3}{2},2}$. 
    The metric, defined in \cite{DohrerReiterSchumacherCompleteRiemannianMetric2025}, is real analytic in $\Mfld$.
\end{corollary}
 
\begin{theorem}\label{thm:Fredholm_Hessian}
	If $\gamma$ is a critical point of the arclength restricted energy $\TPp|_{\AL^{s}}$,
    then the second variation of the energy, restricted to $\T_\gamma \AL^s$, 
     \[
        \D^2 (\TPp)_\gamma \vert_{\T_\gamma\AL^{s}}: \T_\gamma \AL^s\times \T_\gamma \AL^s \rightarrow \R
     \]
     and the Hessian
    \begin{align*}
        &\mathit{Hess}^\AL\TPp_\gamma=\D^2 \left( \TPp \vert_{\AL^s}\right)_{\gamma}:T_\gamma \AL^s \times T_\gamma \AL^s \rightarrow \R,\\
        &\hspace{20pt}(v,w) \mapsto D^2(\TPp)_{\gamma}(v,w)+\D(\TPp)_{\gamma}(h_{12}(v,w)),
    \end{align*}
     induce Fredholm operators with index zero. 
     Here
     $
     \mathit{Hess}^\AL \TPp_\gamma: T_\gamma \AL^s \times T_\gamma \AL^s \rightarrow \R 
     $ denotes the Hessian of $ \TPp $ with respect to $ \AL^s $, and $h_{12} $ denotes the second fundamental form of $\AL^s\subset H^s$. 
\end{theorem}
This leads to the desired \L{}ojasiewicz-Simon gradient inequality.

\begin{theorem}[\L{}ojasiewicz-Simon gradient inequality]\label{thm:LSMaintheorem}
	 Let $ \gamma\in\AL^s $ be a critical point of the restricted energy $ \TPp\vert_{\AL^s} $. 
     Then there are constants $ Z>0 $, $ \delta\in(0,1] $ and $ \theta\in \left[\frac{1}{2},1\right) $, such that for any $ \eta \in \AL^s$ with $ \pdist_{\AL^s}(\gamma,\eta) <\delta$ the following inequality holds true.
    \begin{equation}
    \begin{aligned}
        \left\|\nabla^{\AL}_{\eta}\TPp(\eta)  \right\|_{T_\eta \AL^s}=\left\| \D \left( \TPp\vert_{\AL^s} \right)_{\eta} \right\|_{(T_\eta \AL^s)^*}\ge Z\left|\TPp(\gamma)-\TPp(\eta) \right|^{\theta}.
    \end{aligned}
    \end{equation}
\end{theorem}

As a consequence we can strengthen improve the strong subconvergence of \cite[Theorem 1.7]{FrechesSchumacherSteenebruggevonderMoselPalaisSmaleConditionGeometric2025} and conclude the full convergence of the constrained gradient flow. 

\begin{theorem}\label{thm:strong_convergence_Gradflow}
    Let $ \xi:[0,\infty)\to \AL^s_0 $ be the solution of the Cauchy problem 
    \begin{equation}
    \begin{aligned}\label{eq:CauchyProblemTP}
        \xi'(t)=-\nabla^{\AL^s_{0}}_{\xi(t)} \TPp(\xi(t)) \text{ with } \xi(0)=\gamma_0,
    \end{aligned}
    \end{equation}
    for $ \gamma_0\in\AL^s_{0} $.
    Then $ \xi $ converges strongly in $ \AL_0^s $, for $ t\to\infty $,  to a critical point $ x_\infty $ of $ \TPp\vert_{\AL^s_0} $.
\end{theorem}  

\tocless\subsection{Preliminaries on the arc-length manifold}
In \cite[Chapter 3]{FrechesSchumacherSteenebruggevonderMoselPalaisSmaleConditionGeometric2025} the manifold of arc-length parametrized curves in $H^s(\Domain,\R^n)$ is investigated. 
In the rest of this article, we shall refer to curves parametrized by arc-length by \textit{arc-length curves}.
We give a short review of the results that are of importance for the present paper.

In \cite[Theorem 3.1]{FrechesSchumacherSteenebruggevonderMoselPalaisSmaleConditionGeometric2025} it was shown that $ \AL^s $ is a smooth Riemannian submanifold of $ H^{s}(\Domain,\R^n) $, where the Riemannian structure is given by the standard inner product on $ H^s $.
Adapting the arguments in \cite[Proposition 4.9]{OkabeSchraderConvergenceSobolevGradient2023}, one obtains analyticity of $ \AL^s $. 
For a fixed point $\gamma\in \AL^s$, the tangent space is given by 
\begin{align*}
	\T_\gamma \AL^s=\left\{h\in H^s(\Domain,\R^n)|\langle\gamma'(x),h'(x)\rangle=0 \text{ for all }x\in \Domain \right\}. 
\end{align*}

Another important manifold is given by $ \AL^s_{0}$.
Constraining the problem onto this manifold prevents the curves from traveling wildly through the ambient space along the gradient flow of geometric energies.
By \cite[Theorem 3.2]{FrechesSchumacherSteenebruggevonderMoselPalaisSmaleConditionGeometric2025} this is an analytic Riemannian submanifold of $ \AL^s $.
Its tangent space is given by 
\begin{align*}
	\T_\gamma \AL^s_0=\left\{h\in T_\gamma\AL^s|h(0)=0 \right\}. 
\end{align*}
The manifold $ \AL^s_0 $ was introduced in \cite{FrechesSchumacherSteenebruggevonderMoselPalaisSmaleConditionGeometric2025}, in order to prove the Palais--Smale condition for the generalized tangent point energies, as translations on $ \R^n $ are not compact. 

Clearly, for every $\gamma\in \AL^s$ we have that $\mathcal{L}(\gamma)=\int_{\Domain}\left|\gamma'(\tau) \right|\di \tau=1$, and critical points of $ \TPp\vert_{\AL^s_0} $ are also critical for $\TPp\vert_{\AL^s}$, see \cite[Corollary 3.13]{FrechesSchumacherSteenebruggevonderMoselPalaisSmaleConditionGeometric2025}.

The restriction to arclength curves is a natural assumption, at least for the following two reasons:\\
Firstly, one can show that critical points of the arclength restricted energy are indeed also critical for the length constrained problem (see \cite[Theorem 3.10]{FrechesSchumacherSteenebruggevonderMoselPalaisSmaleConditionGeometric2025}). 
This relies on the invariance under reparametrization and smoothness of critical points, proven in \cite[Theorem 1.8]{FrechesSchumacherSteenebruggevonderMoselPalaisSmaleConditionGeometric2025}.
Thus, by the Lagrange multiplier theorem (see \cite[Corollary 3.5.29]{abrahamManifoldsTensorAnalysis1988}), there exists a constant $\lambda=\lambda(\gamma)\in\R$ such that a critical point satisfies the Euler--Lagrange equation
\begin{align*}
	0=\D\TPp_\gamma v+\lambda\D \mathcal{L}_\gamma(v)\text{ for all } v\in H^s(\Domain,\R^n). 
\end{align*}
Secondly, one can show that the set of injective regular curves, $\Mfld$, deformation retracts onto $\AL^s$, such that the sets are homotopy equivalent, see \cite[Theorem 3.7]{FrechesSchumacherSteenebruggevonderMoselPalaisSmaleConditionGeometric2025}. 
One can even show that this deformation retraction is equivariant, see \cite{freches_2025}.
Thus the two sets are equivariantly homotopy equivalent. 
Note that, without any sort of length constraint, one cannot expect to find critical points at all, as the tangent point energy is negatively homogeneous. 
So, by restricting the problem to the arclength manifold, we lose neither topological nor geometrical nor analytical information.  

\tocless\subsection{Preliminaries on \L{}ojasiewicz--Simon inequality}
The main ingredient in this paper will be a \L{}ojasiewicz--Simon gradient inequality for the restricted energy $ \TPp\vert_{\AL^s} $. 
Strong convergence of the flow will be deduced by standard arguments.
There are plenty similar results to be found in the literature in the context of geometric analysis of curves, especially concerning the Euler--Bernoulli  bending energy under various constraints (for example in \cite{dallacquaOjasiewiczSimonGradientInequality2016}, \cite{ruppExistenceConvergenceLengthpreserving2020b} or \cite{OkabeSchraderConvergenceSobolevGradient2023}).
See also \cite{BlattGradientFlowHara2018} for an application of \L{}ojasiewicz--Simon inequality for the $ L^2 $-gradient flow for O'hara energies. 

While the inequality itself usually appears in more or less the same form, there exist many theorems on when it necessarily holds.
This started with the original work by \L{}ojasiewicz (see \cite{łojasiewicz1965ensembles}), followed by the first generalization to the Hilbert setting due to Simon (see \cite{SimonAsymptoticsforNonlinearEvolutionEquations83}) and the generalizations of Chill (see \cite{chillLojasiewiczSimonGradient2003a}), Huang (see \cite{huangGradientInequalitiesApplications2006}) and Feehan-Maridakies (see \cite{FeehanMaridakisLojasiewiczSimonGradientInequalities2020}).
We employ the version due to Feehan and Maridakies.

\begin{theorem}[{{\cite[Theorem 1]{FeehanMaridakisLojasiewiczSimonGradientInequalities2020}}}]
    \label{thm:LSInequalityFeehan_Maridakies}
    Let $ X\subset X^* $ be a continuous embedding of a Banach space into its dual space. Let $ U\subset X $ be an open subset, let $ \EL:U\to\R $ be an analytic function, and let  $x_\infty $ be a critical point of $ \EL $, such that $ \EL'(x_\infty)=0 $. Assume that $ \EL''(x_\infty):X\to X^* $ is a Fredholm operator of index zero.\\
    Then there exist constants $ Z\in(0,\infty),\sigma\in(0,1] $ and $ \theta\in\left[\frac12.1\right) $ such that the following holds true.
    If $ x\in U $ satisfies 
    \[\left\|x-x_\infty \right\|_{X}<\sigma,\]
    then 
    \[\left\|\EL'(x)   \right\|_{X^*}\ge Z\left|\EL(x)-\EL(x_{\infty}) \right|^{\theta}.\]
\end{theorem}

Note that this result a priori only holds in Banach spaces. However, since the result is only local in nature, one can transfer it to Riemannian submanifolds of Hilbert spaces under certain prerequisites, as in \cite[Proposition 4.13]{OkabeSchraderConvergenceSobolevGradient2023}.
We outline our setting and give a sketch of the proof. The precise and complete proof can be found in \aref{sec:LSineq_section}.

Let $ \HL $ be a Hilbert space, $ \OL\subset\HL $ an open subset, $ \EL:\OL\to\R $ an analytic function, $ \ML\subset\OL $ an analytic Riemannian submanifold of $ \HL $. Let $ \widetilde{\EL}=\EL\vert_{\ML} $ be the restricted energy. 
For a critical point $ x_\infty $ of $ \widetilde{\EL} $ take an analytic chart $ \phi: \UL\to V $ of $ \ML $ containing $ x_\infty $. Then the localized energy $ E: V \to \R, v\mapsto \widetilde{\EL}(\phi^{-1}(v)) $ is analytic and satisfies the suppositions of \cref{thm:LSInequalityFeehan_Maridakies} if and only if $ (\widetilde{\EL})''(x_\infty) $ is Fredholm with index zero. 
This is the key difficulty, as calculating the Hessian of a manifold constrained function is in general far from trivial. Okabe and Schrader resolved this obstruction by considering a modified energy instead of the classical elastic energy, see \cite[Proof of Prop. 4.12]{OkabeSchraderConvergenceSobolevGradient2023}. 
Luckily, we found a convenient workaround in the context of geometric energies on spaces of closed immersed curves based on \cite[Theorem 3.10]{FrechesSchumacherSteenebruggevonderMoselPalaisSmaleConditionGeometric2025} that would also work for other geometric energies on curves, such as the elastic energy, integral Menger curvature or O'Hara energies. 
In fact, the Fredholm property of $ ( \widetilde{\EL} )''(x_\infty) $ can be shown in a way completely analogous to the case of $ \TPp $, but analyticity seems to be much easier for $ \TPp $ compared to $ \mathrm{intM}^{p,2} $ or $ \mathscr{O}^{\alpha,1} $. 
For the elastic energy, analyticity and the Fredholm property has been shown in \cite{OkabeSchraderConvergenceSobolevGradient2023} for a slightly more general manifold. 
In case of the arclength constraint our adapted tools would work as well. In fact, the Fredholm property becomes trivial using our workaround and the regularity results from \cite[Theorem 4.4]{FrechesSchumacherSteenebruggevonderMoselPalaisSmaleConditionGeometric2025}.\\
Note that the choice of $ \HL=H^{s}(\Domain,\R^n) $ and $ \OL=\Mfld $ is very convenient for us, since we are interested in an ODE-type Sobolev gradient flow in the corresponding energy space.

For applications to geometric $ L^2 $-gradient flows our approach might not be suitable. 
Instead, one might consider the refined constrained \L{}ojasiewicz--Simon gradient inequality shown in \cite[Theorem 1.4]{RuppLojasiewiczSimonGradientInequality2020}. 
Note that Rupp's approach only relies on the Fredholm property of $ \EL''(x_\infty) $ which avoids the difficulties with the constrained Hessian $ \left( \EL\vert_{\ML} \right)''(x_\infty) $ in an elegant way. 
However, the suppositions on the operator $ \mathcal{G} $ in Rupp's result make an application in our case unfeasible since our constraints are not finite dimensional.

\section{Analyticity}
In this section we show, that $ \TPp $ is analytic on on $\Mfld$:
\begin{equation}\label{equ:DefEnergy}
	\TPp(\gamma)
		= 
			\int_{\R/ \Z}
			\int_{-1/2}^{1/2}
				\frac{
					\abs{\Delta_{x+w,x} \gamma- D_\gamma\gamma(x)\inner{ D_\gamma\gamma(x), \Delta_{x+w,x} \gamma}}^2
					}
					{\abs{\Delta_{x+w,x} \gamma}^{2s+1}}
					\abs{\gamma'(x)}
					\abs{\gamma'(x+w)}
					\dd w
				\dd x
			.
\end{equation}
In the above, we denote differentiation w.r.t. arclength by $D_\gamma h= \frac{1}{\abs{\gamma'(x)}}h'(x)$ and the difference by $\Delta_{y,x} h= h(y)-h(x)$ .
We will use this notation throughout the rest of the article.
As $ \Mfld $ is an open subset of a Hilbert space, we use the following definition of analyticity. 

\begin{definition}[{\cite{AnalyticFunctionsOnBanachSpacesWhittlesey}}]
Let $(X,\norm{\cdot}_X), (Y,\norm{\cdot}_Y)$ be Banach spaces.
A function 
$f:X \rightarrow Y$ is \textit{(real) analytic around} $x_o\in X$ if there is an open neighborhood $B_R(x_0)$ and $a_n \in L_n(X,Y),c_n\in \R$ such that
\[
	f(x)=\Sigma_{n\in \N} c_n a_n (x-x_0,...., x-x_0) \; \forall x\in B_R(x_0)
		\text{ and }
	\Sigma_{n\in \N} \abs{c_n}\norm{a_n}_{\text{op}} r^n <\infty \;\forall r\in (0,R).
\]
In the above, $L_n(X,Y)$ denotes the space of all n-linear, continuous maps from $X$ to $Y$.
\end{definition}
The following basic results on real analytic functions will be quite useful in the following.

\begin{theorem}[{p. 1079, \cite{AnalyticFunctionsOnBanachSpacesWhittlesey}}]
	Let $ V,W,X $ be Banach spaces, $ D\subset V $ and $ E\subset W $ be open and $ f:D\to W,g:E\to X $ be analytic with $ f(D)\subset E $. Then $ g\circ f $ is analytic. 
\end{theorem}
The most important class of analytic operators is the following. 
\begin{example}
	Let $ \ell\in \N $ and $ V_1,\cdots,V_{\ell},W $ be Banach spaces. 
	If $ a:V_1,\cdots,V_\ell\to W $ is continuous and multilinear, then it is analytic. 
\end{example}
The proof follows directly from the definition of analyticity, see \cite[Example 2.3]{RuppLojasiewiczSimonGradientInequality2020}.

The strategy to show analyticity is, to dissect the integrand of our energy functional into smaller and smaller pieces until we are able to apply the previous results.\\
First note that $\TPp(\gamma)=\norm{f(\gamma)(\cdot)}_{\Lebesgue[2]}^2$. We proof that $f: \Mfld \rightarrow L^2$ is an analytic map.
The claim follows then, from the fact that the map $f\mapsto \norm{f}_{\Lebesgue[2]}^2$ is analytic as a map $\Lebesgue[2]\rightarrow \R$.\\ 
Here, $ f $ is given as 
\begin{align*}
	\TPp(\gamma)
		&=
			\int_\Domain
				\int_{-\frac{1}{2}}^{\frac{1}{2}}
				\left(  
					\frac{
						(
							\Delta_{x+w,x}\gamma- D_\gamma \gamma(x)\inner{D_\gamma\gamma(x), \Delta_{x+w,x}\gamma}
						)
					}{\abs{\Delta_{x+w,x}\gamma}^{s+1/2}}
				\abs{\gamma'(x)}^{1/2}
				\abs{\gamma'(x+w)}^{1/2}
				\right)^2
			\dd (w,x)
		\\
		&=:
			\int_\Domain
				\int_{-\frac{1}{2}}^{\frac{1}{2}}
					\abs{
						f(\gamma)(x,w)
					}^2
				\dd (w,x)
\end{align*}
We split $f$, by adding and subtracting $D_\gamma \gamma(x)\int_{x}^{x+w} \abs{\gamma'(\theta)}\dd \theta$, 
so we can rearrange the numerator in the following way:
\begin{align*}
	&\Delta_{x+w,x}\gamma
		-D_\gamma \gamma(x)\inner{D_\gamma\gamma(x), \Delta_{x+w,x}\gamma}
	\\
	&=
		\Delta_{x+w,x}\gamma
		-
		D_\gamma \gamma(x)\int_{x}^{x+w} \abs{\gamma'(\theta)}\dd \theta
		+
		D_\gamma \gamma(x)\int_{x}^{x+w} \abs{\gamma'(\theta)}\dd \theta
		-
		D_\gamma\gamma(x)\inner{D_\gamma\gamma(x), \Delta_{x+w,x} \gamma}
	\\
	&=
		\Delta_{x+w,x}\gamma
		-
		D_\gamma \gamma(x)\int_{x}^{x+w} \abs{\gamma'(\theta)}\dd \theta
		-
		D_\gamma\gamma(x)
		\inner{
			D_\gamma\gamma(x),
			\Delta_{x+w,x}\gamma
			-
			D_\gamma \gamma(x)
			\int_{x}^{x+w} \abs{\gamma'(\theta)}\dd \theta
		}
\end{align*}
Hence, $f$ can be written as $f(\gamma)=\mathcal{F}(\gamma, \gamma) \Lambda(\gamma)\psi(\gamma)$ with
\begin{align*}
	\mathcal{F}(\gamma, h)
	&=
		\frac{1}{\abs{w}^{s+1/2}}
			\left( 
				\Delta_{x+w,x} h
				- 
				D_\gamma h(x)
				\int_{x}^{x+w} 
					\abs{\gamma'(\theta)}
				\dd \theta
			\right)
	\\
		&
		\hspace{0.5cm}- 
			\frac{1}{\abs{w}^{s+1/2}}
			\left(
				D_\gamma h(x)
				\inner{
					D_\gamma\gamma(x), 
					\Delta_{x+w,x} \gamma
					- 
					D_\gamma\gamma(x)
					\int_{x}^{x+w} 
						\abs{\gamma'(\theta)}
					\dd \theta
				}
			\right),
	\\
	\Lambda(\gamma)
	&=
		(
			\frac{\abs{w}}{\abs{\Delta_{x+w,x} \gamma}}
		)^{(2s+1)/2},
	\\
	\psi(\gamma)
	&=
		(\abs{\gamma'(x+w)}\abs{\gamma'(x)})^{1/2}
		.
\end{align*}
In the above, we interpret $\mathcal{F}$ as a mapping $\mathcal{F}:\Mfld\times \Bessel[s] \rightarrow \Lebesgue[2](\Domain\times \intervaloo{-\frac{1}{2}, \frac{1}{2}})$ and $\Lambda : \Mfld \rightarrow \Lebesgue[\infty](\Domain \times \Domain)$.
Observing that terms repeat, we dissect $\mathcal{F}$ some more.
\[
	\mathcal{F}(\gamma,h)
	= 
	H_1(\gamma,h)
	-H_2(\gamma,h)
\] 
with 
$H_1(\gamma):\Bessel[s](\Domain, \AmbSpace)\rightarrow \Lebesgue[2](\Domain\times (-\frac{1}{2},\frac{1}{2}), \AmbSpace),$
given by 
\begin{align*}
	H_1(\gamma):
		h\mapsto
			\frac{1}{\abs{w}^{s+1/2}}
			\left(
				\Delta_{x+w,x} h
				-
				D_\gamma h(x)
					\int_{x}^{x+w} 
						\abs{\gamma'(\theta)}
					\dd \theta
			\right)
\end{align*}
and 
$H_2:\Mfld\times \Bessel[s] \rightarrow \Lebesgue[2](\Domain\times (-\frac{1}{2},\frac{1}{2}), \AmbSpace)$
by
\begin{align*}
	H_2: (\gamma,h) 
		&\mapsto 
			\frac{1}{\abs{w}^{s+1/2}}
			\left(
				D_\gamma h(x)
				\inner{
					D_\gamma\gamma(x), 
					\Delta_{x+w,x} \gamma
					- 
					D_\gamma\gamma(x) 
					\int_{x}^{x+w}
						\abs{\gamma'(\theta)}
					\dd \theta
				}
			\right)
			\\
		&=
			D_\gamma h(x)
			\inner{D_\gamma\gamma(x), H_1(\gamma,\gamma)}.
\end{align*}
Since $h\mapsto H_1(\gamma)(h)$ and $h\mapsto H_2(\gamma)(h)$ are linear for fixed $\gamma \in \Mfld$, it suffices to show that they are bounded.
Note that we dissected in a way, that the dependencies in $ h $ are linear, and the ones in $ \gamma $ are nonlinear but harmless.
Hence, we divided the problem into showing boundedness of operators and proving that these operators depend analytically on $\gamma$.
The building-blocks are linked through summation and multiplication, which are analytic if they are continuous.
We now show that differentiation w.r.t arclength is analytic.
\begin{lemma}\label{lem:DiffArcAnalytic}
	The map $\gamma \mapsto D_\gamma= \frac{1}{\abs{\gamma'(x)}}\partial_x$
	is analytic as a map $\Mfld\rightarrow L(\Bessel[s], \Bessel[s-1])$.
\end{lemma}
\begin{beweis}
	Let $\gamma_0 \in \Mfld$ be arbitrary and choose $R>0$ such that for all $\gamma \in B_R(\gamma_0)$, the following properties hold true.
	\[
		2 \norm{\gamma_0'}_{\Lebesgue[\infty]}
		\geq 
		\abs{\gamma'(x)}
		\geq 
		\frac{1}{2}\inf_{x\in \Domain}\abs{ \gamma_0'(x)} 
		\text{ for all }x\in \Domain
		\text{  and  }
		\gamma 
		\text{ is an embedding }
		.
	\]
	Since $\partial_x \in L(\Bessel[s], \Bessel[s-1])$ is linear and bounded, it is analytic.
	Furthermore, multiplication $\Bessel[s-1](\Domain,\R)\times\Bessel[s-1](\Domain,\AmbSpace)\to\Bessel[s-1](\Domain,\AmbSpace), (u,v)\mapsto uv$ is bilinear and bounded, since $ s>\frac32 $, thus analytic.
	Thus we need to show that $\gamma\mapsto \frac{1}{\abs{\gamma'}} \in H^{s-1}(\Domain, \R)$ is analytic.
	The smoothness of the map was already shown in \cite[Lemma 3.5]{DohrerReiterSchumacherCompleteRiemannianMetric2025}.
	\\
	We know that $\phi_1:\gamma \mapsto \abs{\gamma'(x)}^2$ is analytic, due to its bilinearity and the fact that $\Bessel[s-1]$ is a Banach-Algebra.
	Furthermore, we know that $\phi_1(\gamma)$ is uniformly bounded from below for all $\gamma \in B_R(\gamma_0)$.
	In addition, $x\mapsto x^{-1/2}$ is analytic on $(c,\infty)$ for all $c>0$.
	Hence, $\gamma\mapsto (\phi_1(\gamma))^{-1/2}$ is analytic.
	This implies that $(\gamma, h)\mapsto D_\gamma h$ is an analytic map $\Mfld\times \Bessel[s] \rightarrow \Bessel[s-1]$.
\end{beweis}\\
The next lemma concerns $\psi$. 
\begin{lemma}\label{lem:VolumeElementAnalytic}
	The map $\psi: \Mfld \rightarrow \Lebesgue[\infty](\Domain \times (\frac{-1}{2}, \frac{1}{2}),\R)$ defined by
	\[
		\psi:\gamma\mapsto (\abs{\gamma'(x)}\abs{\gamma'(x+w)})^{1/2}
	\]
	is analytic.
\end{lemma}
\begin{beweis}
	Let $\gamma_0 \in \Mfld$ be arbitrary.
	Once again, choose $R>0$ such that for all $\gamma \in B_R(\gamma_0)$, it holds that $ \gamma $ is injective and 
	\[
		2 \norm{\gamma_0'}_{\Lebesgue[\infty]}
		\geq 
		\abs{\gamma'(x)}
		\geq 
		\frac{1}{2}\inf_{x\in \Domain}\abs{ \gamma_0'(x)} 
		\text{ for all }x\in \Domain.
	\]
	Such a $ R>0 $ exists by \cite[Lemma 2.76]{steenebruggeSpeedcontrollingGradientFlows2023}.
	Since $s>\frac{3}{2}$, we can use an analogous argument as in the lemma before and conclude that
	\[
		\phi_2:\Bessel[s]\rightarrow L^\infty(\Domain, \R)
		,
		\text{ given by }
		\phi_2:\gamma \mapsto \abs{\gamma'(x)}^{1/2}
	\]
	is analytic.
	Furthermore, the inclusion $\Lebesgue[\infty](\Domain, \R)\rightarrow \Lebesgue[\infty](\Domain \times (\frac{-1}{2}, \frac{1}{2}),\R)$
	is linear and bounded.
	Hence, the map $\phi_2$ is analytic.
	This implies that $\psi$ is the composition and multiplication of analytic maps.
	Therefore, $\psi$ itself is analytic. 
\end{beweis}\\
In the following lemma, we show that $H_1$ depends analytically on $\gamma$ as map $\Mfld \rightarrow L(H^s, L^2)$.
\begin{lemma}\label{lem:H1Analytic}
	The map $H_1:\Mfld \rightarrow  L(H^s, L^2(\Domain \times (\frac{-1}{2}, \frac{1}{2})))$ given by $\gamma \mapsto H_1(\gamma)$ is analytic.
\end{lemma}
\begin{beweis}
	Let $\gamma_0\in \Mfld$ be arbitrary and choose $R>0$ such that
	\[
		2\norm{\gamma_0'}_\infty \geq \abs{\gamma'(x)}\geq \frac{1}{2} \inf \abs{\gamma_0'}
		\; \text{ for all }
		\gamma \in B_R(\gamma)
	.
	\]
	Since $\gamma_0$ was chosen arbitrary, it suffices to show that $\gamma\mapsto H_1(\gamma)$ is analytic in $B_R(\gamma_0)$.
	Let $h\in H^s$ be arbitrary.
	We split $H_1(\gamma)h$ into its building blocks.
	\begin{align*}
		H_1(\gamma)h (x,w)
		&=
			\frac{1}{\abs{w}^{s+1/2}}
			(
				\Delta_{x+w,x} h
				- 
				D_\gamma h(x) 
					\int_x^{x+w} 
						\abs{\gamma'(\theta)}
					\dd \theta
			)
		\\
		&=
			\frac{1}{\abs{w}^{s+1/2}}
			(
				\int_x^{x+w}
					(
						D_\gamma h(\theta)
						-
						D_\gamma h(x)
					)
				\abs{\gamma'(\theta)}
				\dd \theta
			)
		\\
		&=
			\frac{1}{\abs{w}^{s+1/2}}
			(
				\int_x^{x+w}
					(
						D_\gamma h(\theta)
						\abs{\gamma'(\theta)}
						-
						D_\gamma h(x)
						\abs{\gamma'(\theta)}
						\pm 
						D_\gamma h(x)
						\abs{\gamma'(x)}
					)
				\dd \theta
			)
		\\
		&=
			\frac{1}{\abs{w}^{s+1/2}}
			(
				\int_{x}^{x+w}
					(
						h'(x)
						-
						h'(\theta)
					)
				\dd \theta
				+
				D_\gamma h(x)
				\int_{x}^{x+w}
					(
						\abs{\gamma'(x)}
						-
						\abs{\gamma'(\theta)}
					)
				\dd \theta
			)
		\\
		&=
			\phi(h')(x,w)
			-
			D_\gamma h(x)
			\phi( \abs{\gamma'})(x,w)
	\end{align*}
	where $\phi(h)(x,w)=\frac{1}{\abs{w}^{s+1/2}}
			(
				\int_{x}^{x+w}
					(
						h'(x)
						-
						h'(\theta)
					)
				\dd \theta
			)$.
\\
We immediately observe, that 
$
	D_{(\cdot)}: B_R(\gamma_0) \rightarrow L(H^s, H^{s-1}), 
	\gamma \mapsto D_\gamma
$ and 
$
	\abs{\cdot}
		\circ 
	\frac{\dd}{\dd t}
		: 
			B_R(\gamma_0)
			\rightarrow 
			H^{s-1}
$ are analytic.
Furthermore, $\phi$ is linear and does not depend on $\gamma$.
Therefore, it suffices to show that it is a bounded linear map.
The proof uses similar techniques as \cite{Blatt:2013:energyspacestangentpointenergiesa}.
Let $k\in H^{s-1}$ be arbitrary.
\begin{align*}
	\norm{\phi(k)}_{L^2}^2
	&=
		\int_\Domain
			\int_{-1/2}^{1/2}
			\frac{
				\abs{
				\int_{x}^{x+w}
					k(\theta)
					-
					k(x)
				\dd \theta
				}^2
			}
			{\abs{w}^{2s+1}}
		\dd w \dd x
	\\
	&\leq
		\int_\Domain
			\int_{-1/2}^{1/2}
			\frac{
				\int_{0}^{1}
					\abs{
					k(x+ \alpha w)
					-
					k(x)
					}^2
				\dd \alpha
			}
			{\abs{w}^{2s-1}}
		\dd w \dd x
	\\
	&=
		\int_0^1
			\abs{\alpha}^{2s-1}
		\int_\Domain
		\int_{-1/2}^{1/2}
			\frac{
				\abs{k(x+\alpha w)- k(x)}^2
			}
			{\abs{\alpha w}^{2s-1}}
		\dd w \dd x \dd \alpha
	\\
	&=
		\int_0^1
			\abs{\alpha}^{2s-1}
		\int_\Domain
		\int_{-\alpha/2}^{\alpha/2}
			\frac{
				\abs{k(x+z)- k(x)}^2
			}
			{\abs{z}^{2s-1}}
		\dd z \dd x \dd \alpha
	\\
	&\leq
		\int_0^1 \alpha ^{2s-1} \dd \alpha
		\;
		[k]_{s-1}^2
\end{align*}
Therefore $\phi: H^{s-1} \rightarrow L^2(\Domain\times (-\frac{1}{2}, \frac{1}{2}))$ is bounded and linear.
Since $\phi$ does not depend on $\gamma$, we know that $\gamma \mapsto \phi$ is analytic.
Furthermore, since $\frac{\dd }{\dd t}\in L(H^s, H^{s-1})$, we conclude that $\gamma \mapsto \phi \circ \frac{\dd}{\dd t}\in L(H^s, L^2(\Domain\times(\frac{-1}{2}, \frac{1}{2})))$ is analytic.
Additionally, we know that $\gamma \mapsto \phi(\abs{\gamma'})\in L^2$ is analytic.
Combining the above with the fact, that multiplications $M:L^2 \rightarrow L(L^\infty, L^2)$ are bilinear and bounded, we conclude that
\[
	\gamma 
	\mapsto
	D_\gamma (\cdot )
	\frac{
		\int_{x}^{x+w}
			\abs{\gamma'(\theta)} 
		\dd \theta
	}{\abs{w}^{s+1/2}}
	\in 
	L(H^s, L^2(\Domain \times (\frac{-1}{2}, \frac{1}{2})))
\]   
is analytic.
Therefore,
\[
	\gamma \mapsto
	\phi \circ \frac{\dd }{\dd t}
	-
	D_\gamma (\cdot) \phi(\abs{\gamma'})
\]
is analytic as map $\Mfld\rightarrow L(H^s, L^2(\Domain\times(\frac{-1}{2}, \frac{1}{2})))$.
\end{beweis}
\begin{remark}
	As a direct consequence of the above lemma, we conclude that the map $\gamma \mapsto H_1(\gamma)\gamma$ is analytic.
	Since the reverse is not true, we prove the claim.
	Since $H_1$ is analytic in $\gamma$, we know that there exists $R>0$ and $a_n \in L_n(H^s, L(H^s,L^2))$ such that
	\[
		H_1(\gamma)
		=
		\sum_{n\in \N} 
			a_n(\gamma- \gamma_0)^n
		\text{ for all }
		\gamma \in B_R(\gamma_0)
	\]
	with $\sum_n \norm{a_n}_{\text{op}} r^n <\infty$ for all $0<r<R$.
	Hence
	\[
		H_1(\gamma)\gamma = \sum_n a_n (\gamma-\gamma_0)^n (\gamma)= \sum_n a_n (\gamma-\gamma_0)^n (\gamma_0) + a_n (\gamma-\gamma_0)^{n}(\gamma-\gamma_0)
	\]
	We define $b_0\ceq a_0(\gamma_0)$ and for $n\geq 1$ $b_n\ceq a_n(\gamma-\gamma_0)^n(\gamma_0)+ a_{n-1}(\gamma-\gamma_0)^{n-1}(\gamma-\gamma_0)$.
	We are left with showing the absolute convergence of the series.
	Since this is a local property, we may assume that $0<R<1$.
	\begin{align*}
		\sum_{n\in \N} \norm{b_n}_{\text{op}} r^n
		&=
			\norm{b_0}_{\text{op}}
			+
			\sum_{n\geq 1} \norm{a_n (\cdot)(\gamma_0) + a_{n-1}}_{\text{op}} r^n
		\\
		&\leq
			\norm{a_0(\gamma_0)}_{\text{op}}
			+
			\sum_{n\geq 1} \norm{\gamma_0}\norm{a_n}_{\text{op}} r^n + \norm{a_{n-1}}_{\text{op}}r^n
		\\
		&\leq
			\norm{a_0(\gamma_0)}_{\text{op}}
			+
			\norm{\gamma_0}\sum_{n\geq 1} \norm{a_n}_{\text{op}} r^n + \sum_{n\geq 1} \norm{a_{n-1}}_{\text{op}}r^{n-1}
		<\infty
	\end{align*}
\end{remark}
In the following lemma, we make use of the fact, that $H_2$ can be written as the product of analytic functions and $H_1$. 
Note that the multiplication operator, $M:\Lebesgue[\infty] \times \Lebesgue[2]\rightarrow \Lebesgue[2]$, is bounded and bilinear, thus analytic.

\begin{lemma}
	The function $H_2:\Mfld \times \Bessel[s] \rightarrow \Lebesgue[2](\Domain\times (\frac{-1}{2},\frac{1}{2}), \AmbSpace)$ is analytic.
\end{lemma}	
\begin{beweis}
	Let $(\gamma,h)\in  \Mfld\times \Bessel[s] $ be arbitrary.
	$H_2$ can be written as
	\[
		H_2(\gamma,h)
			= 
				M
				\left(
					D_\gamma h(x), 
					M(D_\gamma \gamma(x), H_1(\gamma)\gamma)
				\right)
		.
	\]
	The map $(\gamma,h)\mapsto D_\gamma h$ is analytic, due to \aref{lem:DiffArcAnalytic}.
	Furthermore, the embeddings $i_1: \Bessel[s-1] \hookrightarrow \Lebesgue[\infty](\Domain), i_2 :\Lebesgue[\infty](\Domain) \hookrightarrow \Lebesgue[\infty](\Domain\times (\frac{-1}{2},\frac{1}{2}), \AmbSpace)$ are linear and continuous, hence analytic.
	Using \aref{lem:H1Analytic}, we conclude that $\gamma\mapsto H_1(\gamma)\gamma$ is analytic.
	Since compositions of analytic maps are analytic, we now conclude that $(h,\gamma) \mapsto H_2(h,\gamma)$ is analytic. 
\end{beweis}

Lastly, we investigate $\Lambda$.

\begin{lemma}
	The function $\Lambda: \Mfld \rightarrow \Lebesgue[\infty]( \Domain \times (\frac{-1}{2}, \frac{1}{2}))$, defined via 
	\[
		\gamma 
		\mapsto 
		\left(
			\frac{\abs{w}}
			{
				\abs{ \gamma(x+w)-  \gamma(x)}
			}
		\right)^{\frac{2s+1}{2}}
	\]
	is analytic in $\Mfld$.
\end{lemma}
\begin{beweis}
	Let $\gamma_0 \in \Bessel[s]$ be an arbitrary, immersed embedding 
	and as before $ R>0 $ such that every $ \gamma\in B_R(\gamma_0) $ is injective and 
	\[
		2\norm{\gamma_0'}_{\Lebesgue[\infty]}
		> 
		\abs{ \gamma'(x)} 
		> 
		\frac{1}{2} 
		\inf_{x\in \Domain} \abs{\gamma_0'(x)} 
		\coloneqq c(\gamma_0).
	\]
	Therefore, all $\gamma \in B_R(\gamma_0)$ are embeddings with 
	\[
		\varphi(\gamma)(x,w)
		=
		\frac{ 
			\abs{ \gamma(x+w)-\gamma(x)}^2}
			{\abs{w}^2}
			\in (c(\gamma_0)^2, 4\norm{\gamma_0'}_{\Lebesgue[\infty]}^2).
	\]
	Since $\varphi$ is bilinear and bounded in $B_R(\gamma_0)$, we know that it is analytic.
	Now, we make use of the fact, that $h: y \mapsto y^{\frac{-(2s+1)}{4}}$ is an analytic and bounded function from $(c,\infty) \rightarrow \R$ for all $c>0$.
	This yields that
	\[
		h\circ \varphi : 
		B_R(\gamma_0) \cap \AL^s 
		\rightarrow \Lebesgue[\infty](\Domain \times (\frac{-1}{2}, \frac{1}{2}))
	\]
	is analytic.
\end{beweis}
\\
In order to prove the main theorem of this section, it remains to show that the manifold of arclength curves itself is analytic.
\begin{lemma}\label{lem:ALsAnalyticHigherS}
    Let $s> \frac32$.
    The manifold $\AL^s\subset \Bessel[s](\Domain, \AmbSpace)$ is an analytic submanifold.
\end{lemma}
\begin{beweis}  
    The claim follows the same lines, as \cite[Theorem 3.1]{FrechesSchumacherSteenebruggevonderMoselPalaisSmaleConditionGeometric2025}, where one can always use the embedding $ H^s\hookrightarrow C^1 $, as $ s>\frac{3}{2} $. 
    The defining submersion is given by $ \Sigma:\Mfld\to H^{s-1}(\Domain,\R), \gamma\mapsto\ln(\left|\gamma'\right|) $. 
    The operator $ \gamma\to\gamma' $ is analytic. 
    Furthermore, $ \ln $ and the euclidean norm are analytic away from zero, thus $ \Sigma $ is analytic. 
    By the preimage theorem and the analytic version of the implicit function theorem, the manifold is analytic. 
    See \cite[Proposition 4.9]{OkabeSchraderConvergenceSobolevGradient2023} for a similar result in the case $ s=2 $.
\end{beweis}    

We now collect the fruits of our hard work and proof the first main result, namely that $\TPp$ is analytic (see \cref{thm:analyticityTPp}).
\\
\begin{beweis}[ of \cref{thm:analyticityTPp}]
	We write
	\[
		\TPp(\gamma)
		=
			\norm{
				M
				\left(
					M
					\pars{
						H_1(\gamma)(\gamma)- H_2(\gamma, \gamma)
					,
						\Lambda(\gamma)
					}
				, \psi(\gamma)
				\right)
				}_{\Lebesgue[2]}^2
			\quad,
	\]
	whereas before, $M$ denotes the multiplication operator $L^2\times L^\infty\rightarrow L^2, (u,v)\mapsto uv$.
	The sum of analytic functions is analytic, therefore $H_1(\cdot)(\cdot)-H_2(\cdot,\cdot):\Mfld\times \Bessel[s] \rightarrow \Lebesgue[2](\R/\Z\times (\frac{-1}{2},\frac{1}{2}), \AmbSpace)$ is analytic.
	As $ M $ is analytic, we have that
	\[
		(\gamma, h) \mapsto 
		M
		\left(
		M
				\pars{
					H_1(\gamma)\gamma- H_2(h, \gamma)
					,
					\psi(\gamma)
					}
		,
		\Lambda(\gamma)
		\right)
	\]
	is analytic, as a map $\Mfld\times \Bessel[s] \rightarrow \Lebesgue[2](\R/\Z\times (\frac{-1}{2},\frac{1}{2}), \AmbSpace)$.
	Lastly we use, that the map $ f\to \left\|f \right\|^2_{L^2} $ is analytic, since 
	\[
		\norm{ f}_{L^2}^2 = \inner{f, f}_{L^2}
	\]
	and the inner product is a continuous, bilinear map. This finishes the proof.
\end{beweis}\\
We deduce that the restriction of the tangent-point energy to $\AL^s$ is analytic.
\begin{corollary}
	The map $\TPp\vert_{\AL^s}: \AL^s \rightarrow \R,\; p=2s+1$ is analytic.
\end{corollary}
\begin{beweis}
	By the previous result, we know that $\TPp$ is real analytic on $\Mfld$. Furthermore, $\AL^s$ is an analytic submanifold of $\Mfld$.
	Since the restriction of analytic maps to analytic submanifolds is analytic, we conclude the proof.
\end{beweis}\\
As a consequence, we can state the following corollary, which extends the theory developed by the first author in \cite{DohrerReiterSchumacherCompleteRiemannianMetric2025},
and allows us to use tools from the theory of analytic functions in order to analyze functions on $(\Mfld,G)$.
Furthermore, this implies that the second fundamental form of analytic submanifolds like $\AL^s$ w.r.t $(\Mfld,G) $ is well-defined.
\begin{corollary}
	The Riemannian metric $G$ on $\Mfld$ defined in \cite{DohrerReiterSchumacherCompleteRiemannianMetric2025} is analytic.
\end{corollary}
\begin{beweis}
	Recall that the metric is given by
	\begin{align*}
		G_\gamma(h,k)
			&= \inner{h,k}_{\Lebesgue[2](\Domain, \AmbSpace, \abs{\gamma'}\dd x)}
				+
				\inner{D_\gamma h,D_\gamma k}_{\Lebesgue[2]( \AmbSpace, \abs{\gamma'}\dd x)}
			\\
			\quad 
			&+
				\iint_{\Domain^2}
					\inner{\Op{\gamma}{s}h, \Op{\gamma}{s}k}
					\frac{\abs{\gamma'(x)}\abs{\gamma'(y)}}{\abs{\gamma(y)-\gamma(x)}}
				\dd (y,x)
			\\
			\quad 
			&+
				\iint_{\Domain^2}
					\abs{\Op{\gamma}{s}\gamma}^2
					\frac{\inner{h(y)-h(x), k(y)-k(x)}}{\abs{\gamma(y)-\gamma(x)^2}}
					\frac{\abs{\gamma'(x)}\abs{\gamma'(y)}}{\abs{\gamma(y)-\gamma(x)}}
				\dd (y,x)
			\\
			\quad 
			&+
				\iint_{\Domain^2}
					\abs{\Op{\gamma}{s}\gamma}^2
					(
						\inner{D_\gamma h(x), D_\gamma k(x)}
						+
						\inner{D_\gamma h(y), D_\gamma k(y)}
					)
				\frac{\abs{\gamma'(x)}\abs{\gamma'(y)}}{\abs{\gamma(y)-\gamma(x)}}
				\dd (y,x)
				,
	\end{align*}
where $\Op{\gamma}{s}h= \frac{1}{\abs{\Delta_{y,x} \gamma}^s}(\Delta_{y,x} h- D_\gamma h(x)\inner{D_\gamma\gamma(x), \Delta_{y,x} \gamma})$.
We did already show, that the occurring geometric $H^1$ metric is analytic in $\gamma$.
Carefully writing out the terms, we realize that
\[
	\frac{\Op{\gamma}{s}h}{\abs{\gamma(y)-\gamma(x)}^{1/2}}\abs{\gamma'(x)}^{1/2}\abs{\gamma'(y)}^{1/2}
	=
		(H_1(\gamma, h)- H_2(\gamma, h))\psi(\gamma)\Lambda(\gamma).
\]
We did already show that the above expression is analytic in $\gamma$ and $h$.
Furthermore, multiplying the above by analytic $\Lebesgue[\infty]$ functions is analytic.
Therefore, we are left with showing that the following functions are analytic in $\gamma, h,k$:
\begin{itemize}
	\item $F_1:(\gamma, h,k)\mapsto \frac{\inner{h(y)-h(x), k(y)-k(x)}}{\abs{\gamma(y)-\gamma(x)}^2}$
	\item $F_2:(\gamma, h,k)\mapsto (\inner{D_\gamma h(x), D_\gamma k(x)}+ \inner{D_\gamma h(y), D_\gamma k(y)})$
\end{itemize}
We have already shown that $F_1,F_2$ are bounded and linear in $h,k$ as maps $\Bessel[s]\rightarrow \Lebesgue[2](\Domain\times\Domain)$.
Furthermore, we know that $\gamma \mapsto D_\gamma \in L(\Bessel[s], \Bessel[s-1])\subset L(\Bessel[s], \Lebesgue[2](\Domain\times\Domain))$ is analytic.
Therefore, it suffices to show that $F_1$ is analytic in $\gamma$.
Using that
\[
	F_1(\gamma, h, k)
	=
		\frac{
			\inner{h(y)-h(x), k(y)-k(x)}
		}
		{\dist[\Domain](y,x)^2}
		(
			\frac{\dist[\Domain](y,x)}
			{\abs{\gamma(y)-\gamma(x)}}
		)^2,
\]
 we can directly conclude that $F_1$ is analytic in $\gamma$. 
\end{beweis}
\section{The Fredholm property}\label{sec:Fredholm_Geometric}

In this section we prove, that the Hessian of the generalized tangent point energy with respect to $ \AL^s $ induces a Fredholm operator of index zero. 
When restricting to the arclength manifold, the second derivative of the energy at a fixed point $ \gamma $ is a symmetric bilinear form
\begin{align*}
	&\mathit{Hess}^{\AL^s}\TPp_\gamma(v,w)=\D^2 \left( \TPp \vert_{\AL^s}\right)_{\gamma}: \T_\gamma \AL^s\times \T_\gamma \AL^s\to \R,\\
	&\hspace{20pt} (v,w)\mapsto D^2(\TPp)_{\gamma}(v,w)+\D(\TPp)_{\gamma}(h_{12}(v,w)),
\end{align*}
where $ h_{12}: \T_\gamma\AL^s\times \T_\gamma\AL^s\to \left( \T_\gamma\AL^s \right)^\perp $ denotes the second fundamental form of $ \AL^s $, see \cite[XIV,\S1 and \S2]{LangFundamentalsDifferentialGeometry1999}.  One cannot hope, that the second term vanishes, as $ \gamma $ is only critical on $ \AL^s $. On the other hand, the second fundamental form is not easy, to directly calculate. However we can use a little workaround: 
As mentioned earlier, for every critical point $ \gamma $ of $\TPp|_{\AL^s}$, there exists $ \lambda\in \R $ such that $\gamma$ is critical for the functional
\begin{align*}
	\mathcal{EL}=\TPp(\cdot)+\lambda \mathcal{L}(\cdot).
\end{align*}
Note that for arclength constrained curves $f$ the two energies $\mathcal{EL}$ and $\TPp$ only differ by a constant. Since $\gamma$ is critical for $\mathcal{EL}$, the second derivative, $\D^2(\mathcal{EL}|_{\AL^s})$ equals the restriction of $\D^2(\mathcal{EL})$ to $(\T_\gamma\AL^s)^2$, as $ \D(\mathcal{EL})_{\gamma}(h_{12}(v,w))=0 $.\\
On the other hand we have by the linearity of $ \mathit{Hess}^{\AL^s} $, that
\begin{align*}
	\mathit{Hess}^{\AL^s}(\mathcal{EL})_\gamma=\mathit{Hess}^{\AL^s}(\TPp)_\gamma+\mathit{Hess}^{\AL^s}(\mathcal{L})_\gamma,
\end{align*} 
where the second term on the right side vanishes, as the length functional is constant on $\AL^s$.\\
In summary, we can calculate $ \D^2(\TPp_{\gamma}\vert_{\AL^s}) = \D^2(\mathcal{EL})_\gamma\vert_{(\T_\gamma\AL^s)^2}$, which is still technical, but not nearly as bad, as it could have been.\\
We will make excessive use of \cite[Theorem 4.15]{FrechesSchumacherSteenebruggevonderMoselPalaisSmaleConditionGeometric2025},
where the second author and his collaborators proved the $C^\infty$-smoothness of critical points of $\TPp\vert_{\AL^s}$.
This is crucial, as it allows us, to use fractional integration by parts, in order to show compactness of lower order terms, by "throwing more differentiation over" onto critical points $\gamma$.
For a more compact exposition of the following arguments, we adapt the notation from \cite{DohrerReiterSchumacherCompleteRiemannianMetric2025}.
We now proceed to show, that the Hessian differs from their strong Riemannian metric only by compact perturbations, and is thus Fredholm.
The metric is closely related to $\TPp$ and $D\TPp$. 
Let $s\in \intervaloo{\frac{3}{2},2},\gamma \in \Mfld, h,k \in T_\gamma \Mfld$.
We have
\begin{align*}
	G_\gamma(h,k)&= \inner{h,k}_{L^2(\Domain, \abs{\gamma'}\dd x)}
					+
					\inner{D_\gamma h,D_\gamma k}_{L^2(\Domain, \abs{\gamma'}\dd x)}
					+
					\sum_{i=1}^3 B_\gamma^i(h,k)
\end{align*}
where
\begin{align*}
	B_\gamma^1 (h,k)
	&=
		\iint_{\Domain^2}
			\inner{\Op{\gamma}{s}h, \Op{\gamma}{s}k}
		\dd \mu_\gamma
	\\
	B_\gamma^2(h,k)
	&=
		\iint_{\Domain^2}
			\abs{\Op{\gamma}{s}\gamma}^2
			\frac{\inner{h(y)-h(x),k(y)-k(x)}}{\abs{\gamma(y)-\gamma(x)}^2}
		\dd \mu_\gamma
	\\
	B_\gamma ^3(h,k)
	&=
		\iint_{\Domain^2}
			\abs{\Op{\gamma}{s}\gamma}^2
			(
				\inner{D_\gamma h(x), D_\gamma k(x)}
				+
				\inner{D_\gamma h(y), D_\gamma k(y)}
			)
		\dd \mu _\gamma
\end{align*}
with
\[
\Op{\gamma}{s} h
	=
		\frac{L_\gamma h }{\abs{\gamma(y)-\gamma(x)}^s}   
	=
	\frac{1}{\abs{\gamma(y)-\gamma(x)}^s}
	(h(y)-h(x)-D_\gamma h(x)\inner{D_\gamma \gamma(x), \gamma(y)-\gamma(x)})
	,
\] 
$\dd \Omega_\gamma(x,y)= \abs{\gamma'(x)}\abs{\gamma'(y)} \dd (x,y)$ and $\dd\mu_\gamma= \frac{\dd\Omega_\gamma}{\abs{\gamma(y)-\gamma(x)}}$.
A straight forward computation yields that
\[
	\TPp(\gamma)=
	B_\gamma^1(\gamma, \gamma)
	\text{ for }
	p=2s+1.
\]
We differentiate once in direction $h\in T_\gamma \AL^s$ and obtain the following
\begin{align*}
	D\TPp(\gamma) h
	&=
		2B_\gamma ^1(\gamma, h)
		-
		p B_\gamma^2 (\gamma, h)
		+
		B_\gamma^3(\gamma, h)
	.
\end{align*}

Furthermore, it was already shown that 
\begin{align*}
	D L \vert_\gamma(\varphi)h
	&=
		-D_\gamma h(x)(
			\inner{D_\gamma\gamma(x), L_\gamma\varphi}
			+
			\inner{D_\gamma\varphi(x),L_\gamma\gamma}
		).
\end{align*}
These identities simplify differentiating the energy.
We write $(DB^i(\gamma) \phi)(u,v)$ as short form of $D(\gamma \mapsto B_\gamma ^i(u,v))(\gamma) \phi$ in order to indicate that we only differentiate the dependency on the base point.
We now differentiate once more, this time in direction $k\in T_\gamma \AL^s$.
Here, we use, that we can write the first derivative as the sum of three $\gamma$-dependent bilinear forms, evaluated at $\gamma$ and $h$.
Hence, we use the product rule and obtain the following.
\begin{align*}
	D^2\TPp(\gamma)(h,k)
	&=
		2B_\gamma^1 (h,k)
		-
		p B_\gamma^2(h,k)
		+
		B_\gamma^3(h,k)
	\\
	&+
		2 (D B^1(\gamma)k)(\gamma, h)
		-
		p
		(DB^2(\gamma)k)(\gamma, h)
		+
		(DB^3(\gamma)k)(\gamma,h)
\end{align*}
Since the Hessian is symmetric, we know that if it is Fredholm, it has index 0.
We are going to show, that $B_\gamma^1$ induces a Fredhom-operator.
The remaining bilinear forms induce compact operators $T_\gamma\AL^s \rightarrow (T_\gamma \AL^s)^*$.
This suffices, since compact perturbations of Fredholm-operators are Fredholm, see \cite[Chapter XVII Corollary 2.6]{LangRealFunctionalAnalysis1993}.
We now start by sorting out the compact terms.
\begin{lemma}\label{lem:B2compact}
	Let $\gamma \in \AL^s$ be arbitrary.
	The map $(h,k)\mapsto B_\gamma^2(h,k)$ is a compact, bilinear map $\Bessel[s]\rightarrow \R$.
\end{lemma}
\begin{beweis}
	First, we observe that $\gamma\in \AL^s$ implies that $\TPp(\gamma)<\infty$ (see \cite[Proposition 2.4]{blattRegularityTheoryTangentpoint2015}).
	Combining this with \cite[Proposition 2.3]{blattRegularityTheoryTangentpoint2015}, we conclude that $\BiLip(\gamma)<\infty$.
	This is a straight forward computation.
	\begin{align*}
		\abs{B_\gamma^2(h,k)}
		&=
		\abs{
			\iint_{\Domain^2}
				\abs{ \Op{\gamma}{s} \gamma}^2
				\frac{ \inner{\Delta_{y,x} h, \Delta_{y,x} k}}{\abs{\Delta_{y,x} \gamma}^2}
				\dd \mu_\gamma
		}
		\\
		&\leq
			\BiLip(\gamma)^2
			\norm{ D_\gamma h}_{\Lebesgue[\infty]}	
			\norm{ D_\gamma k}_{\Lebesgue[\infty]}
			\TPp(\gamma)
		\\
		&=
			\Bilip(\gamma)^2
			\norm{ h'}_{\Lebesgue[\infty]}	
			\norm{ k'}_{\Lebesgue[\infty]}
			\TPp(\gamma)
	\end{align*}
	Since the embedding $\Bessel[s] \hookrightarrow C^1$ is compact, there is nothing left to show.
\end{beweis}\\
In the next Lemma, we have to make use of the additional regularity of $\gamma$, which is granted by $\gamma$ being a critical point.
\begin{lemma}
	Let $\gamma\in \AL^s$ be a critical point of $\TPp$.
	The map $(h,k)\mapsto (DB ^2 (\gamma)k)(\gamma,h)$ is a bilinear, compact map $\Bessel[s] \times \Bessel[s]\rightarrow \R$.
\end{lemma}
\begin{beweis}
	First, we differentiate $\gamma \mapsto B_\gamma^2(\psi_1,\psi_2)$ in direction $k$ for $\psi_1,\psi_2 \in \Bessel[s]$.
	For $u\in \AmbSpace\setminus \{0\}$, we denote by $P_u^\perp (v)=v- \frac{\inner{u,v}}{\abs{u}^2}$ the projection onto the orthogonal complement of $\R u$.
	\begin{align*}
		&(D B^2(\gamma) k)(\psi_1,\psi_2)
		\\
		&=
			\frac{\dd}{\dd \epsilon}\bigg\vert_{\epsilon=0}
			\iint_{\Domain^2}
				\frac{
					\abs{
						P^{\perp}_{(\gamma+\epsilon k)'(x)}(\Delta_{y,x} \gamma+\epsilon k)
					}^2
				}{\abs{ \Delta_{y,x} \gamma+ \epsilon k}^{2s+3}}
				\inner{ \Delta_{y,x} \psi_1, \Delta_{y,x} \psi_2}
			\dd \Omega_{\gamma+ \epsilon k}
		\\
		&=
			-(2s+3)
			\iint_{\Domain^2} \abs{\Op{\gamma}{s} \gamma}^2
				\frac{\inner{\Delta_{y,x} \psi_1, \Delta_{y,x} \psi_2}}{\abs{\Delta_{y,x} \gamma}^2}
				\frac{\inner{\Delta_{y,x} \gamma, \Delta_{y,x} k}}{\abs{\Delta_{y,x} \gamma}^2}
			\dd \mu_\gamma
		\\
		&
		\quad
			+2
			\iint_{\Domain^2} 
				\inner{\Op{\gamma}{s}\gamma, \Op{\gamma}{s} k}
				\frac{
					\inner{\Delta_{y,x} \psi_1, \Delta_{y,x} \psi_2}
					}
				{\abs{\Delta_{y,x} \gamma}^2}
			\dd \mu_\gamma
			\\
		&
		\quad+
			\iint_{\Domain^2}
				\abs{\Op{\gamma}{s}\gamma}^2 
				\frac{\inner{\Delta_{y,x} \psi_1, \Delta_{y,x} \psi_2}}{\abs{\Delta_{y,x} \gamma}^2}
				\left(
					\inner{D_\gamma \gamma (x), D_\gamma k(x)}
					+
					\inner{D_\gamma\gamma(y), D_\gamma k(y)}
				\right)
				\dd \mu_\gamma
	\end{align*}
The third summand vanishes, since $k\in T_\gamma \AL^s$.
We now analyze the remaining terms.
Inserting $\psi_1=\gamma, \psi_2= h$ and choosing $\epsilon= \frac{2-s}{2}$, one obtains
\begin{align*}	
	 (D B^2 (\gamma)k)(\gamma,h)
	 &=
	 	-(2s+3)
		\iint_{\Domain^2}
			\abs{\Op{\gamma}{s}\gamma}^2
			\frac{\inner{\Delta_{y,x} \gamma, \Delta_{y,x} h}}{\abs{\Delta_{y,x} \gamma}^2}
			\frac{\inner{\Delta_{y,x} \gamma, \Delta_{y,x} k}}{\abs{\Delta_{y,x} \gamma}^2}
		\dd \mu_\gamma
		\\
	&\quad
	+2
		\iint_{\Domain} 
			\inner{ \Op{\gamma}{s} \gamma, \Op{\gamma}{s} k}
			\frac{
					\inner{\Delta_{y,x} \gamma, \Delta_{y,x} h}
					}
				{
					\abs{\Delta_{y,x} \gamma}^2
				}
		\dd \mu_\gamma
	\\
	&=
		-(2s+3)
		\iint_{\Domain^2}
			\abs{\Op{\gamma}{s}\gamma}^2
			\frac{\inner{\Delta_{y,x} \gamma, \Delta_{y,x} h}}{\abs{\Delta_{y,x} \gamma}^2}
			\frac{\inner{\Delta_{y,x} \gamma, \Delta_{y,x} k}}{\abs{\Delta_{y,x} \gamma}^2}
		\dd \mu_\gamma
		\\
	&\quad
	+2
		\iint_{\Domain} 
			\inner{ \Op{\gamma}{s+\epsilon} \gamma, \Op{\gamma}{s-\epsilon} k}
			\frac{
					\inner{\Delta_{y,x} \gamma, \Delta_{y,x} h}
					}
				{
					\abs{\Delta_{y,x} \gamma}^2
				}
		\dd \mu_\gamma
	\quad 
	.
\end{align*}
With a similar argument as for $B_\gamma^2 (h,k)$, we bound the first term by $C(\gamma,s)\norm{h'}_{\Lebesgue[\infty]}\norm{k'}_{\Lebesgue[\infty]}$.
Hence, we conclude that the first term is compact in $h$ and $k$.
In order to conclude that the second term is compact, we estimate it in the following way.
But first, we use again that $\gamma\in \AL^s$ implies that $\TPp(\gamma)<\infty$ and therefore that $\BiLip(\gamma)<\infty$.
\begin{align*}
	&\abs{
		\iint_{\Domain} 
			\inner{ \Op{\gamma}{s+\epsilon} \gamma, \Op{\gamma}{s-\epsilon} k}
			\frac{
					\inner{\Delta_{y,x} \gamma, \Delta_{y,x} h}
					}
				{
					\abs{\Delta_{y,x} \gamma}^2
				}
		\dd \mu_\gamma
	}
	\\
	&
	\quad
	\leq
	\BiLip(\gamma)
	\norm{h'}_{\Lebesgue[\infty]}
		\iint_{\Domain^2}
			\abs{\Op{\gamma}{s+\epsilon}\gamma}
			\abs{\Op{\gamma}{s-\epsilon} k}
		\dd \mu_\gamma
	\\
	&
	\quad
	\leq
		\BiLip(\gamma)
		\norm{h'}_{\Lebesgue[\infty]}
		\sqrt{ \iint_{\Domain^2} 
			\abs{\Op{\gamma}{s+\epsilon} \gamma}^2
			\dd \mu_\gamma}
		\sqrt{ 
			{ \iint_{\Domain^2} 
			\abs{\Op{\gamma}{s-\epsilon} k}^2
			\dd \mu_\gamma}
			}
	\\
	&
	\quad
	\leq
		C(\gamma,s, \epsilon)
		\norm{h'}_{\Lebesgue[\infty]}
		\norm{\gamma}_{\Bessel[s+\epsilon]}
		\norm{k}_{\Bessel[s-\epsilon]}
\end{align*}
In the above, we used \cite[Theorem 4.1]{DohrerReiterSchumacherCompleteRiemannianMetric2025} with $s_1=s-\epsilon, s_2=s+\epsilon$.
Due to \cite[Theorem 4.15]{FrechesSchumacherSteenebruggevonderMoselPalaisSmaleConditionGeometric2025},
we know that $\norm{\gamma}_{\Bessel[{s_2}]}$ is finite.
Combining this bound with the fact, that $\Bessel[s]\hookrightarrow\Bessel[s_1]= \Bessel[s-\epsilon]$ is compact, we conclude that $(h,k) \mapsto (DB ^2 (\gamma)k)(\gamma,h)$ is compact.
\end{beweis}\\

Next, we investigate $B_\gamma^3$ and its derivative.
The occuring terms are either compact or vanish due to $h,k$ being tangential to $\AL^s$.

\begin{lemma}\label{lem: B3compact}
	Let $\gamma \in \Mfld$. The map $(h,k)\mapsto B_\gamma^3(h,k)$ is a bilinear, compact map
	$\Bessel[s]\times \Bessel[s] \rightarrow \R$.
\end{lemma}
\begin{beweis}
	Let $\gamma\in \Mfld$ and $ h,k \in \Bessel[s]$ be arbitrary.
	Due to \cite[Proposition 2.4]{blattRegularityTheoryTangentpoint2015}, we know that $\TPp(\gamma)<\infty$.
	Furthermore, we know that there is a $C>0$, such that
	\[
		C \geq \abs{\gamma'(x)} \geq \frac{1}{C} \text{ for all }
		x\in \Domain
		.
	\]
	Therefore, we can estimate $B_\gamma ^3$ in the following way.
	\begin{align*}
		B_\gamma^3(h,k)
		&=
			\int_{\Domain}
				\int_{-\frac{1}{2}}^\frac{1}{2}
					\abs{
						\Op{\gamma}{s} \gamma
					}^2
					(\inner{D_\gamma h(x), D_\gamma k(x)}+ \inner{D_\gamma h(x+w), D_\gamma k(x+w)})
				\dd \mu_\gamma
		\\
		&\leq
			2C^2 
			\norm{h'}_{\Lebesgue[\infty]}
			\norm{k'}_{\Lebesgue[\infty]}
			\int_{\Domain}
				\int_{-\frac{1}{2}}^\frac{1}{2}
					\abs{
						\Op{\gamma}{s} \gamma
					}^2
					\dd \mu_\gamma
		=
			C \TPp(\gamma)
			\norm{h'}_{\Lebesgue[\infty]}
			\norm{k'}_{\Lebesgue[\infty]}
		\end{align*}
This shows, that $B_\gamma^3: W^{1,\infty}\times W^{1,\infty} \rightarrow \R$ is a bilinear, continuous map.
Hence, $B_\gamma^3: \Bessel[s] \times \Bessel[s]\rightarrow\R$ is bilinear and compact.
\end{beweis}\\

\begin{lemma}\label{lem: DB3compact}
	Let $\gamma \in \AL^s$ and $h,k \in T_\gamma \AL^s$.
	Then
	\[
		(DB^3 (\gamma)k)(\gamma,h)=0
	.
	\]
\end{lemma}
\begin{beweis}
	Let $\gamma \in \AL^s$ and $h,k \in T_\gamma \AL^s$ be arbitrary.
	This implies, that $\inner{\gamma'(x), h'(x)}=\inner{\gamma'(x), k'(x)}=0$ for all $x\in \Domain$.
	Recall that 
	\begin{equation}
	\begin{aligned}\label{eq:B3Formel}
		B_\gamma ^3(\psi_1,\psi_2)
		=
			\int_{\Domain}
				\int_{-\frac{1}{2}}^\frac{1}{2}
					\abs{\Op{\gamma}{s} \gamma}^2
					(
						\inner{D_\gamma \psi_1(x), D_\gamma \psi_2(x)}
						+
						\inner{D_\gamma \psi_1(x+w), D_\gamma \psi_2(x+w)}
					)
				\dd \mu_\gamma.
	\end{aligned}
	\end{equation}

	Differentiating \cref{eq:B3Formel} after $ \gamma $ yields
	\begin{align*}
		&(DB^3 (\gamma)k)(\psi_1,\psi_2)=\\
		&= 
		\int_{\Domain}
				\int_{-\frac{1}{2}}^\frac{1}{2} 
					D\left( 
						\left|\Op{\gamma}{s} \gamma \right|^2
						\frac{\abs{\gamma'(x)} \abs{\gamma'(x+w)}}
						{\abs{\gamma(x+w)-\gamma(x) }} 
					\right)
					(
						\inner{D_\gamma \psi_1, D_\gamma \psi_2}\vert_{x}
						+
						\inner{D_\gamma \psi_1, D_\gamma \psi_2}\vert_{x+w}
					)\di w \di x,
		\\
		&\hspace{0.5cm}
		-2
			\int_{\Domain}
				\int_{-\frac{1}{2}}^\frac{1}{2} 
					\abs{\Op{\gamma}{s}\gamma}^2
					(
						\inner{D_\gamma \psi_1, D_\gamma \psi_2}
						\inner{D_\gamma \gamma, D_\gamma k}
						\vert_{x}
						+
						\inner{D_\gamma \psi_1, D_\gamma \psi_2}
						\inner{D_\gamma \gamma, D_\gamma k}\vert_{x+w}
					)
				\dd \mu_\gamma
	\end{align*}
	and replacing $ \psi_1,\psi_2 $ with $ \gamma,h $ for $ h\in \T_\gamma\AL^s $ yields $ (DB^3 (\gamma)k)(\gamma,h) =0$, as $ \left\langle \gamma'(x),h'(x) \right\rangle=0  $ for all $ x\in \Domain $. 
\end{beweis}\\

We now analyze the derivative of $B_\gamma^1$. 
Using the same techniques as established in the previous proofs and \cite[Theorem 4.15]{FrechesSchumacherSteenebruggevonderMoselPalaisSmaleConditionGeometric2025}, we observe that the occurring terms are compact.

\begin{lemma}
	Let $\gamma\in \AL^s$ be a critical point of $\TPp, p=2s+1$.\\
	The map $(h,k) \mapsto (DB^1 (\gamma)k)(\gamma,h)$ is a bilinear, compact map $\Bessel[s]\times \Bessel[s]\rightarrow\R$.
\end{lemma}
\begin{beweis}
	First, we observe that $\gamma\in \AL^s$ implies that $\TPp(\gamma)<\infty$ (see \cite[Proposition 2.4]{blattRegularityTheoryTangentpoint2015}).
	Combining this with \cite[Proposition 2.3]{blattRegularityTheoryTangentpoint2015}, we conclude that $\BiLip(\gamma)<\infty$.
	We start by writing out the derivative itself. 
	\begin{align*}
		(DB^1 (\gamma)k)(\psi_1,\psi_2)
		&=
			-(2s+1)
				\iint_{\Domain^2}
					\inner{\Op{\gamma}{s}\psi_1,\Op{\gamma}{s}\psi_2}
					\frac{
						\inner{\Delta_{y,x}k,\Delta_{y,x}\gamma}
					}{\abs{\Delta_{y,x}\gamma}^2}
				\dd \mu_\gamma
		\\
		&-
			\iint_{\Domain^2}
				\frac{1}{\abs{\Delta_{y,x} \gamma}^{2s}}
				(
					-\inner{F_\gamma \psi_1, D_\gamma \psi_2}
					(
						\inner{D_\gamma k, F_\gamma\gamma}
						+
						\inner{D_\gamma\gamma, F_\gamma k}
					)
				)
				\dd \mu_\gamma
		\\
		&-
		\iint_{\Domain^2}
		\frac{1}{\abs{\Delta_{y,x} \gamma}^{2s}}
		(
			-\inner{F_\gamma \psi_2, D_\gamma \psi_1}
			(
				\inner{D_\gamma k, F_\gamma\gamma}
				+
				\inner{D_\gamma\gamma, F_\gamma k}
			)
		)
		\dd \mu_\gamma
	\end{align*}
Inserting $\psi_1=\gamma, \psi_2=h$, we obtain the following expression.
\begin{align*}
	(DB^1 (\gamma)k)(\gamma,h)
	&=
		-(2s+1)
			\iint_{\Domain^2}
				\inner{\Op{\gamma}{s}\gamma,\Op{\gamma}{s}h}
				\frac{
					\inner{\Delta_{y,x}k,\Delta_{y,x}\gamma}
				}{\abs{\Delta_{y,x}\gamma}^2}
			\dd \mu_\gamma
	\\
	&-
		\iint_{\Domain^2}
			\frac{1}{\abs{\Delta_{y,x} \gamma}^{2s}}
			(
				-\inner{F_\gamma \gamma, D_\gamma h}
				(
					\inner{D_\gamma k, F_\gamma\gamma}
					+
					\inner{D_\gamma\gamma, F_\gamma k}
				)
			)
			\dd \mu_\gamma
	\\
	&-
	\iint_{\Domain^2}
	\frac{1}{\abs{\Delta_{y,x} \gamma}^{2s}}
	(
		-\inner{F_\gamma h, D_\gamma \gamma}
		(
			\inner{D_\gamma k, F_\gamma\gamma}
			+
			\inner{D_\gamma\gamma, F_\gamma k}
		)
	)
	\dd \mu_\gamma
\end{align*}
We investigate the occurring summands one by one, starting with the first one.
We choose $\epsilon= \frac{2-s}{2}$ and compute
\begin{align*}
	\abs{
		\iint_{\Domain^2}
			\inner{\Op{\gamma}{s}\gamma,\Op{\gamma}{s}h}
			\frac{
				\inner{\Delta_{y,x}k,\Delta_{y,x}\gamma}
			}{\abs{\Delta_{y,x}\gamma}^2}
		\dd \mu_\gamma
	}
&\leq
	\BiLip(\gamma)
	\norm{D_\gamma k}_{\Lebesgue[\infty]}
		\iint_{\Domain^2}
			\abs{\Op{\gamma}{s+\epsilon}\gamma}
			\abs{\Op{\gamma}{s-\epsilon}h}
		\dd \mu_\gamma
\\
&\leq
	C(\gamma)
	\norm{k'}_{\Lebesgue[\infty]}
	\norm{h}_{\Bessel[s-\epsilon]}
	\norm{\gamma}_{\Bessel[s+\epsilon]}
	.
\end{align*}
This is compact in $h,k$, since $\Bessel[s]\hookrightarrow\Bessel[s-\epsilon]$ 
and $\Bessel[s]\hookrightarrow C^1$ are compact embeddings.
Once again, thanks to \cite[Theorem 4.15]{FrechesSchumacherSteenebruggevonderMoselPalaisSmaleConditionGeometric2025},
we know that $\norm{\gamma}_{\Bessel[s+\epsilon]}$ is finite.
Next, we deal with the terms involving $F_\gamma \gamma$.
The trick is the same as before: One can shift some "differentiation" from h (respectively k)
to $\gamma$ through moving some of the singularity from $\Op{\gamma}{s}h$ over to $\Op{\gamma}{s}\gamma$.
Nonetheless, we perform one of the two estimates.
\begin{align*}
	&\abs{
		\iint_{\Domain^2}
			\frac{1}{\abs{\Delta_{y,x}\gamma}^{2s}}
			\inner{
				F_\gamma\gamma,
				D_\gamma h
			}
			\inner{
				D_\gamma \gamma,
				F_\gamma k
			}
			\dd \mu_\gamma
	}
	\leq
		\norm{D_\gamma h}_{\Lebesgue[\infty]}
		\iint_{\Domain^2}
			\abs{\Op{\gamma}{s+\epsilon}\gamma}
			\abs{\Op{\gamma}{s-\epsilon}k}
		\dd \mu_\gamma
	\\
	&\leq
		C(\gamma)
		\norm{h'}_{\Lebesgue[\infty]}
		\norm{k}_{\Bessel[s-\epsilon]}
		\norm{\gamma}_{\Bessel[s+\epsilon]}
\end{align*}
With a similar argument as before, this is compact in $h$ and $k$.
The summand involving $\inner{F_\gamma h, D_\gamma\gamma}\inner{D_\gamma k, F_\gamma\gamma}$
can be dealt with analogously.\\
The term involving $\inner{F_\gamma\gamma, D_\gamma h}\inner{D_\gamma k, F_\gamma\gamma}$ can be bounded
in terms of the $\Lebesgue[\infty]$-norms of $h$ and $k$:
\begin{align*}
	\abs{
		\iint_{\Domain} 
			\inner{\Op{\gamma}{s}\gamma, D_\gamma h}
			\inner{\Op{\gamma}{s}\gamma, D_\gamma k}
		\dd \mu_\gamma
	}
	&\leq
		\norm{D_\gamma h}_{\Lebesgue[\infty]}
		\norm{D_\gamma k}_{\Lebesgue[\infty]}
		\iint_{\Domain^2}
			\abs{\Op{\gamma}{s}\gamma}^2
		\dd \mu_\gamma
	\\
	&\leq
		\norm{h'}_{\Lebesgue[\infty]}
		\norm{k'}_{\Lebesgue[\infty]}
		\TPp(\gamma).
\end{align*}
The above is compact in $h$ and $k$, due $\Bessel[s]\hookrightarrow C^1$ being compact.
Hence, we are left with
\begin{align*}
	\iint_{\Domain^2}
		\inner{\Op{\gamma}{s}h, D_\gamma\gamma}
		\inner{D_\gamma\gamma, \Op{\gamma}{s}k}
	\dd \mu_\gamma
	\quad.
\end{align*}
We claim that this term, although it is less obvious, is compact as well.
First, we use the Hölder-inequality and bound the term by 
\begin{align*} 
	\iint_{\Domain^2}
		\inner{\Op{\gamma}{s}h, D_\gamma\gamma}
		\inner{D_\gamma\gamma, \Op{\gamma}{s}k}
	\dd \mu_\gamma
	&\leq
	(\iint_{\Domain^2}\abs{P^\top_{\gamma'(x)}(\Op{\gamma}{s}h)}^2 \dd \mu_\gamma)^{1/2}
	(\iint_{\Domain^2}\abs{P^\top_{\gamma'(x)}(\Op{\gamma}{s}k)}^2 \dd \mu_\gamma)^{1/2}
.
\end{align*}
In the above, $P^\top_{\gamma'(x)}(v)= \gamma'(x)\inner{\gamma'(x),v}$ denotes the tangential projection along the curve $ \gamma $.
Due to the symmetry in $h$ and $k$, we only investigate one of them.
\begin{align*}
	\inner{L_\gamma h, \gamma'(x)}
	&=
	\inner{\Delta_{x+w,x}h- h'(x)\inner{\gamma'(x), \Delta_{x+w,x}\gamma}, \gamma'(x)}
	\\
	&=
		\inner{\Delta_{x+w,x} h, \gamma'(x)}
	=
		w
		\int_0^1 
			\inner{h'(x+\theta w), \gamma'(x)}
		\dd \theta
	\\
	&=
	w
	\int_0^1 
		\inner{h'(x+\theta w), \gamma'(x)}
		-
		\inner{h'(x+\theta w), \gamma'(x+\theta w)}
	\dd \theta
	\\
	&=
	w 
	\int_0^1 
		\inner{h'(x+\theta w), \Delta_{x+\theta w, x}\gamma'}
	\dd \theta
\end{align*}
Using this, we can estimate our last term in the following manner.
\begin{align*}
	&\iint_{\Domain^2}
			\abs{
				P^T_{\gamma'(x)}(\Op{\gamma}{s}h)
				}^2
			\dd \mu_\gamma
	\leq
		\BiLip(\gamma)^{2s+1}
		\iint_{\Domain^2}
			\frac{
				\abs{
				P^T_{\gamma'(x)}(L_{\gamma}h)}^2
			}{\dist[\Domain][y][x]^{2s+1}}
			\dd (y,x)
	\\
	&=
		\BiLip(\gamma)^{2s+1}
			\int_\Domain
				\int_{-1/2}^{1/2}
					\frac{w^2}{\abs{w}^{2s+1}}
					\abs{
					\int_0^1
						\inner{h'(x+\theta w), \Delta_{x+\theta w,x}\gamma'}
					\dd \theta
					}^2
				\dd w
			\dd x
	\\
	&\leq
		\BiLip(\gamma)^{2s+1}
		\int_\Domain
		\int_{-1/2}^{1/2}
		\int_0^1
			\frac{1}{\abs{w}^{2s-1}}
			\norm{h'}_{\Lebesgue[\infty]}^2
			\abs{\Delta_{x+\theta w, x}\gamma'}^2
		\dd\theta
		\dd w   
		\dd x
	\\
	&\leq
		\BiLip(\gamma)^{2s+1}
		\norm{h'}_{\Lebesgue[\infty]}^2
		\int_0^1 \theta ^{2s-1} \dd \theta
		\int_\Domain
		\int_{-1/2}^{1/2}
			\frac{
				\abs{
				\Delta_{x+z,x}\gamma'}^2
				}{\abs{z}^{2s-1}}
			\dd z
			\dd x
	= C(\theta,\gamma)
	\norm{h'}_{\Lebesgue[\infty]}^2
	\seminorm{\gamma'}_{\Bessel[s-1]}^2
\end{align*}
Once again, since $\Bessel[s]\hookrightarrow C^1$ is compact, the above term is compact in $h$.
\end{beweis}\\
Now that we have managed to estimate all terms, we can prove the Fredholm property of the restricted second variation.

\begin{theorem}\label{thm:D2TPpFredholm}
	Let $\gamma \in \AL^s$ be a critical point of $\TPp\vert_{\AL^s}$.
	Then the second derivative of $\TPp$ constrained to $T_\gamma \AL^s$
	\[
		\left( D^2 \TPp (\gamma) \right)\vert_{T_\gamma\AL^s}: T_\gamma\AL^s \rightarrow (T_\gamma\AL^s )^*
	\]
	is a Fredholm operator of index 0.
	\end{theorem}
\begin{beweis}
	Let $\gamma\in \AL^s$ be a critical point of $\TPp\vert_{\AL^s}$.
	We have already derived the expression for the second derivative.
	Due to the observations before, we know that
	\begin{align*}
D^2\TPp(\gamma)(h,k)
	&=
		2B_\gamma^1 (h,k)
		-
		p B_\gamma^2(h,k)
		+
		B_\gamma^3(h,k)
	\\
	&
	\hspace{0.5cm}+
		2 (D B^1(\gamma)k)(\gamma, h)
		-
		2 p
		(DB^2(\gamma)k)(\gamma, h)
		-
		(DB^3(\gamma)k)(\gamma,h)
	\\
		&= 
			B_\gamma^1(h,k) + K(\gamma)(h,k),
	\end{align*}
	 where $K(\gamma):\Bessel[s]\times \Bessel[s] \rightarrow \R$ is compact.
	Hence, it suffices to show that $B_\gamma ^1$ induces a Fredholm-operator.
	Since
	\[
		G_\gamma(\cdot, \cdot)=
		\inner{\cdot, \cdot}_{\Lebesgue[2](\Domain, \abs{\gamma'} \dd x)}
		+
		\inner{D_\gamma\cdot,D_\gamma \cdot}_{\Lebesgue[2](\Domain, \abs{\gamma'} \dd x)}
		+ \sum_1^3 B_\gamma^i(\cdot,\cdot)
	\]
	is a strong Riemannian metric, it induces a Fredholm-operator on $ \T_{\gamma}\AL^s $.
	We now show, that $G$ differs from $B_\gamma^1$ by a compact perturbation.
	It is easy to see, that 
	$
		\inner{\cdot, \cdot}_{\Lebesgue[2](\Domain, \abs{\gamma'} \dd x)}
	$
		and 
	$
		\inner{D_\gamma\cdot,D_\gamma \cdot}_{\Lebesgue[2](\Domain, \abs{\gamma'} \dd x)}
	$ are compact perturbations, because $\Bessel[s]\hookrightarrow C^1$ is compact. 
	We have proved in \aref{lem:B2compact} and \aref{lem: B3compact} that $B_\gamma ^2$ and $B_\gamma^3$ are compact perturbations.
	This now implies, that
	\[
		D^2\TPp(\gamma)= G_\gamma +\tilde{K}(\gamma)
	\]
	where 
	\[
		\tilde{K}(\gamma)
			=
			K(\gamma) 
			- 
			\inner{\cdot,\cdot}_{\Lebesgue[2](\Domain, \abs{\gamma'} \dd x)}
			- 
			\inner{D_\gamma \cdot,D_\gamma \cdot}_{\Lebesgue[2](\Domain, \abs{\gamma'} \dd x)}
			-B_\gamma ^2
			-B_\gamma^3
	\] is a compact $\Bessel[s]\times \Bessel[s] \rightarrow \R$.
	Since inner-products induce Fredholm-operators of index 0, we conclude that $D^2\TPp(\gamma)$ is a Fredholm-operator of index 0.
\end{beweis}

Unfortunately, this is not sufficient to show a \L{}ojasiewicz-Simon gradient inequality.
In order to do so, we need to prove that the Hessian with respect to the connection on the manifold $\AL^s$, by which we mean $\mathit{Hess}^{\AL^s} \TPp_\gamma=\D^2\left( \TPp\vert_{\AL^s} \right)_{\gamma} $, induces a Fredholm operator.

\begin{beweis}[of \cref{thm:Fredholm_Hessian}]
	Let $\gamma\in \AL^s$ be a critical point of $\TPp \vert_{\AL^s}$.
	Due to our observation at the beginning of the section and our choice of $\lambda$, we know that $\gamma$ is a critical point for
	$\mathcal{EL}: \Mfld\rightarrow \R$.
	Since $\Mfld \subset \Bessel[s]$ is an open subset, we know that the Hessian is just the second variation.
	Furthermore, we already observed that
	\[
		D^2(\TPp\vert_{\AL^s})_\gamma
		=
		D^2( \mathcal{EL})_\gamma \vert_{T_\gamma \AL^s}
		=
		D^2( \TPp )_\gamma \vert_{T_\gamma \AL^s}
		+
		\lambda
		D^2 (\mathcal{L}) _\gamma \vert_{T_\gamma \AL^s}
		.
	\]
	In \aref{thm:D2TPpFredholm}, we have shown, that $D^2( \TPp )_\gamma \vert_{T_\gamma \AL^s}$ induces a Fredholm operator.
	It remains to prove $D^2 (\mathcal{L})_ \gamma: \T_\gamma\AL^s \times \T_\gamma\AL^s \rightarrow \R$ is compact.
	Let $\gamma \in \AL^s$ and $h,k\in \Bessel[s]$ be arbitrary, then one easily computes
	\begin{align*}
		D\mathcal{L}_\gamma (h)
			&=
				\int_{\Domain}
					\inner{D_\gamma \gamma, D_\gamma h}
				\abs{\gamma'(\theta)}
				\dd \theta,
	\end{align*}
	and 
	\begin{align*}
		D^2 (\mathcal{L})_\gamma (h,k)
			&=
				\int_{\Domain}
					\inner{D_\gamma h,D_\gamma k}
				\abs{\gamma'(\theta)}
				\dd \theta
				-
				\int_{\Domain}
					\inner{D_\gamma\gamma, D_\gamma h}
					\inner{D_\gamma \gamma, D_\gamma k}
				\abs{\gamma'(\theta)}
				\dd \theta
	.
	\end{align*}
	For $ \gamma\in \AL^s $ and $ h,k\in \T_\gamma\AL^s $ the above reduces to
	\begin{align}
		\label{eq:length_compact}
		D^2 (\mathcal{L})_\gamma (h,k)
			&=\int_{\Domain}\left\langle h',k' \right\rangle \dd \theta,
	\end{align}
	which can be bound by $C(\gamma) \left\|h \right\|_{C^1}\left\|k \right\|_{C^1} $ and thus, is clearly compact. The claim now follows from \cref{thm:D2TPpFredholm}
\end{beweis}

\section{\L{}ojasiewicz-Simon gradient inequality}\label{sec:LSineq_section}
In this section we prove a \L{}ojasiewicz-Simon gradient inequality for geometric energies on curves in the arclength manifold $ \AL^s $ and apply this result to $ \TPp\vert_{\AL^s} $ in order to prove \cref{thm:LSMaintheorem,thm:strong_convergence_Gradflow}.

\begin{theorem}
    \label{thm:LSineq_abstract}
    Let $ s> \frac{3}{2}$, $ \EL:\Mfld\to\R $ be analytic and invariant under reparametrization. 
    Let $ \gamma_\infty\in\AL^s\cap H^{s}(\Domain,\R^n) $ be a critical point of the constrained energy $ \widetilde{\EL}:=\EL\vert_{\AL^s} $. 
    Assume that $ D^2\EL(\gamma_{\infty})\vert_{\T_{\gamma_\infty}\AL^s}:\T_{\gamma_\infty}\AL^s \to\left( \T_{\gamma_\infty}\AL^s \right)^*$ is a Fredholm operator with index zero. 
    Then there are constants $ Z>0 $, $ \delta\in\intervaloc{0,1} $ and $ \theta\in \left[ \frac{1}{2},1\right) $, such that for $ \eta \in \AL^s$ with $ \pdist_{\AL^s}(\gamma_{\infty},\eta) <\delta$ 
    \begin{equation}
    \begin{aligned}
        \left\|\nabla^{\AL^s}_{\eta}\EL(\eta)  \right\|_{T_\eta \AL^s}=\left\| \D \left( \EL\vert_{\AL^s} \right)_{\eta} \right\|_{(T_\eta \AL^s)*}\ge Z\left|\EL(\gamma_{\infty})-\EL(\eta) \right|^{\theta}.
    \end{aligned}
    \end{equation} 
\end{theorem}
\begin{beweis}
    Since $ \gamma_{\infty}\in H^{s} $ and $ \EL $ is invariant under reparametrisation, the suppositions of \cite[Theorem 3.10]{FrechesSchumacherSteenebruggevonderMoselPalaisSmaleConditionGeometric2025} hold, and there is a $ \lambda\in\R $ such that $ \gamma_{\infty} $ is also critical for $ \EL(\cdot)-\lambda\mathcal{L}(\cdot) $. With the arguments from the beginning of \cref{sec:Fredholm_Geometric} we have
    \begin{equation}
    \begin{aligned}
        D^2
            \left( 
                \EL\vert_{\AL^s} 
            \right)(\gamma_{\infty})
        =
            \left( 
                    D^2\EL(\gamma_{\infty}) 
                    -
                    \lambda 
                    D^2\mathcal{L}(\gamma_{\infty})
                \right)
                \big\vert_{
                    \T_{
                        \gamma_{\infty}
                        }\AL^s
                    },
    \end{aligned}
    \end{equation}
    where the first summand is Fredholm with index zero by assumption and the second is compact by \aref{eq:length_compact}. Thus, $ D^2\left( \EL\vert_{\AL^s} \right)(\gamma_{\infty}) $ is Fredholm with index zero as well.\\
    Now let $ U\subset \AL^s $ be an open neighborhood of $ \gamma_{\infty} $ and choose $ \delta $ small enough, such that $ \eta\in U $, and a chart $ \phi:U\to V $, where $ V $ is a subspace of $ H^s(\Domain,\R^n) $. 
    Define the localized Energy $ E:=\widetilde{\EL}\circ\phi^{-1}:\phi(U)\to \R $. 
    Then $ \D E_v(w)=\D \widetilde{\EL}_{\phi^{-1}(v)}\circ \D\phi^{-1}_v(w) $ and $ \D^2 E_{\phi(\gamma)} = \D^2 \widetilde{\EL}_{\gamma}\left( \D\phi^{-1}_{_{\phi(\gamma)}}(\cdot),\D\phi^{-1}_{_{\phi(\gamma)}}(\cdot)  \right)  $, where we used that $ \gamma $ is a critical point of $ \widetilde{\EL} $.\\
    Since $ \AL^s $ is analytic, $ \phi $ is also analytic and $ \D\phi: \T_\gamma\AL^s\to V $ is an isomorphism of Hilbert spaces.
    Hence, $ E $ is analytic, and its second derivative is Fredholm with index zero. 
    Thus, the energy $ E $ satisfies the suppositions of \cref{thm:LSInequalityFeehan_Maridakies} and there are constants $ \tilde{Z}>0 $, $ \tilde{\delta}\in(0,1] $ and $ \theta\in \left[\frac{1}{2},1\right) $ such that 
    \begin{align*}
        \left\|\D E_{\phi(\eta)} \right\|_{V^*}\ge \tilde{Z}\left|E(\phi(\eta))-E(\phi(\gamma)) \right|^{\theta}=\tilde{Z}\left|\EL(\eta)-\EL(\gamma) \right|^{\theta}.
    \end{align*}
    By assumption $ \D \phi^{-1} $ is continuous.
    Hence, there is a constant $ c>0 $ such that for any $ \eta $ with $ \pdist_{\AL^s}(\gamma,\eta)\le \delta $ and any $ v\in T_{\eta}\AL^s $
    \begin{align*}
        \left\|v \right\|_{T_{\eta}\AL^s}=\left\|v \right\|_{H^{s}}\le c \left\|D\phi_{\eta}v  \right\|_{V}
        .
    \end{align*}
    Therefore, 
    \begin{align*}
        \left\|\D E_{\phi_{\eta}} \right\|_{V*}=\underset{\D \phi_{\eta}(v)\in V}{\sup} \frac{\left|\D E_{\phi(\eta)} \D \phi_{\eta}(v)\right|}{\left\|\D\phi_{\eta}(v) \right\|_{V}}\le \underset{v\in\T_{\eta}\AL^s}{\sup}\frac{\left|\D \widetilde{\EL}_{\eta}v \right|}{c^{-1}\left\|v \right\|_{T_{\eta}\AL^s}}=c\left\|\D\left( \EL\vert_{\AL^s} \right) \right\|_{(\T_\eta \AL^s)^*}.
    \end{align*}
    The existence of $ \delta>0 $ such that if $ \pdist_{\AL^s} (\gamma,\eta)<\delta$ then $ \left|\phi(\eta)-\phi(\gamma) \right|_{V} <\tilde{\delta}$ follows from the continuity of $ \phi $ and because $ \AL^s $ is a submanifold of $ H^{s}\left( \Circle, \R^n \right) $. 
\end{beweis}

\begin{remark}\label[remark]{rem:LSineq_fixedpoint}
    The above theorem also holds if one chooses $ \gamma_{\infty}\in \AL^s_{0} $ and $ \widetilde{\EL}=\EL\vert_{\AL^s_{0}} $ and adjusting the other assumptions to $ \AL^s_{0} $. The proof follows from the exact same steps, where one can still apply \cite[Theorem 3.10]{FrechesSchumacherSteenebruggevonderMoselPalaisSmaleConditionGeometric2025} by virtue of \cite[Corollary 3.13]{FrechesSchumacherSteenebruggevonderMoselPalaisSmaleConditionGeometric2025}. 
\end{remark}
\begin{beweis}[ of \cref{thm:LSMaintheorem}]
    It was shown in \cite{FrechesSchumacherSteenebruggevonderMoselPalaisSmaleConditionGeometric2025}, that critical points of $ \TPp\vert_{\AL^s} $ are smooth. Analyticity of $ \TPp $ was shown in \cref{thm:analyticityTPp} and the assumptions on the restricted second variation were verified in \cref{thm:Fredholm_Hessian} and the claim follows from \cref{thm:LSineq_abstract}. 
\end{beweis}

We want to apply the previous theorem to the gradient flow of the restricted energy $\widetilde{\EL}= \EL\vert_{\AL^s} $ 
\begin{equation}
\begin{aligned}
    \label{eq:Cauchyproblemsolution}
    \xi'(t)=-\nabla^{\AL^s} \EL(\xi(t)) \text{ with } \xi(0)=\gamma_0.
\end{aligned}
\end{equation}
The following result is quite well known, but we carry out the proof anyway for the sake of readability and completeness. To that extend we follow the proof of \cite[Theorem 5.5]{OkabeSchraderConvergenceSobolevGradient2023} very closely, but note that the proof of Okabe and Schrader is more involved as they prove a stronger statement than we do.

\begin{lemma}
	\label{convergence_gradflow}
    Let $ \xi:[0,\infty) \to \AL^s$ be a solution of \aref{eq:Cauchyproblemsolution}, and $ \gamma_{\infty} $ a limit point of $ \xi $. Assume in addition, that the suppositions of \cref{thm:LSineq_abstract} hold. Then $ \xi(t)\to \gamma_\infty $ strongly, as $ t\to \infty $ and $ \gamma_\infty $ is a critical point of $ \EL \vert_{\AL^s}$. 
\end{lemma}
\begin{beweis}
    Choose a sequence $t_k\to\infty$ monotonically, such that $\pdist(\xi(t_k),\gamma_\infty)<\frac{\delta}{2}$ for all $k\in\N$.
    For each $k\in\N$, we define the time $T_k$ as the biggest time, such that $\pdist(\xi(t),\gamma_{\infty})\le\delta$ for $t\in[t_k,T_k)$.
    We define the function 
    \begin{align*}
        H(t)\ceq\left( \EL(\xi(t))-\EL(\gamma_{\infty}) \right)^{1-\theta}
        .
    \end{align*}
    Then $H(t)$ is monotonically decreasing and bounded from below by $0$. Differentiation with respect to $ t $ together with \cref{thm:LSineq_abstract} yield
    \begin{align*}
        -H'(t)  &=-(1-\theta)\left( \EL(\xi(t))-\EL(\gamma_{\infty}) \right)^{-\theta}\partial_{t}\EL(\xi(t))\\
                &=(1-\theta)\left( \EL(\xi(t))-\EL(\gamma_{\infty}) \right)^{-\theta} \left\|\nabla^{\AL^s} \EL(\xi(t)) \right\|^2_{\xi(t)}\\
                &\ge (1-\theta)K\left\|\nabla^{\AL^s}  \EL(\xi(t)) \right\|_{\xi(t)}.
    \end{align*}
    Now integration over $[t_i,T_i]$ yields 
    \begin{align*}
        (1-\theta)K\int_{t_i}^{T_i}\left\|\nabla^{\AL^s}  \EL(\xi(r))  \right\|_{\xi(r)}\di r\le H(t_i)-H(T_i).
    \end{align*}
    Fix any $j$ and set $I:=\bigcup_{i\ge j}[t_i,T_i)$.
    We have
    \begin{align*}
        \int_{I}\left\|\nabla^{\AL^s}  \EL(\xi(r))  \right\|_{\xi(r)}\di r\le \frac{H(t_j)}{(1-\theta)K}, 
    \end{align*}
    since $H(t)$ is strictly decreasing. It remains to show, that there exists a $N\in\N$, such that $T_N=\infty$. \\
    Suppose not. 
    Then for every $i\in\N$ there exists a finite $\T_i$, such that $\pdist(\xi(T_i),p)=\delta$.
    Passing to a subsequence if necessary, we may assume that the intervals $[t_i,T_i)$ are disjoint. 
    We then obtain 
    \begin{align*}
        \delta  &= \pdist(\xi(T_i),\gamma_\infty)\le \pdist(\xi(t_i),\gamma_{\infty})+\pdist(\xi(t_i),\xi(T_i))\\
                &\le  \frac{\delta}{2} + \pdist(\xi(t_i),\xi(T_i))\\
                &\le  \frac{\delta}{2} + \int_{t_i}^{T_i}\left\|\xi'(r) \right\|_{\xi(r)}\di r\\
                &\le  \frac{\delta}{2} + \int_{t_i}^{T_i}\left\|\nabla^{\AL^s}  \EL(\xi(r)) \right\|_{\xi(r)}\di r
    \end{align*}
    However, this implies that $\int_{t_i}^{T_i}\left\|\nabla^{\AL^s}  \EL(\xi(r)) \right\|_{\xi(r)}\di r\ge\frac{\delta}{2}$, which in turn means 
    \begin{align*}
        \int_{I}\left\|\nabla^{\AL^s}  \EL(\xi(r)) \right\|_{\xi(r)}\di r\ge \sum_{i\ge N}^{\infty}\int_{t_i}^{T_i}\left\|\nabla^{\AL^s} \EL(\xi(r)) \right\|_{\xi(r)}\di r\ge \sum_{i\ge N}\frac{\delta}{2}=\infty,
    \end{align*}
    a contradiction. 
    This implies 
    \begin{align*}
        \int_{t_N}^{\infty}\left\|\xi'(r) \right\|\di r<\infty.
    \end{align*}
    The strong convergence of $\xi(t)$ for $t\to\infty$ follows by the auxiliary result below. 
    This concludes the proof.
\end{beweis}
\begin{remark}\label[remark]{rem:convergence_gradflow_fixedpoint}
    Similarly to \cref{rem:LSineq_fixedpoint}, one can prove the previous result for the manifold $ \AL^s_{0} $.
\end{remark}
In \cite{FrechesSchumacherSteenebruggevonderMoselPalaisSmaleConditionGeometric2025} the following condition was introduced
\begin{enumerate}
    \item[\textrm (L**)]\label{condL**} 
	every Cauchy sequence 
	$(x_k)_k\subset\ML$ with $\sup_{k\in\N}\EL(x_k)<\infty$ 
	converges to some 
	$x_\infty\in\ML$ as $k\to\infty$. (\textsc {limit in $\ML$})
\end{enumerate}
With this, we can formulate the following.
\begin{lemma}\label[lemma]{lem:finite_length_trajectories}
   Let $ \xi:(a,b)\to \AL^s $ a curve of finite length and let $ \widetilde{\EL} $ be as above, satisfying \textrm (L**).  
    If $ \left|\EL(\xi(t)) \right|\le C $ for a fixed constant $ C\in\R $, then the limit $ \lim\limits_{t\to b}\xi(t) $ exists in $ \AL^s $. 
\end{lemma}
\begin{beweis}
    \label[lemma]{fintelength_trajectory}
    The method of proof is essentially a reiteration of \cite[Lemma 3.3]{OkabeSchraderConvergenceSobolevGradient2023} with the sole difference, that the manifold $ \AL^s $ is not complete.
    Thus, we need the uniform energy bound in combination with \textrm (L**) in order to avoid self intersections in the limit and obtain convergence.\\
    Take any monotonely increasing sequence $t_k\to b$. 
    One can easily show, that $\xi(t_k)$ forms a Cauchy sequence in $(\AL^s,\pdist)$. 
    Since the energy satisfies the condition \textrm (L**) by assumption we conclude convergence to some $\xi(b)$. 
    Furthermore, the limit is unique, since we can take any other sequence $\tilde{t}_k$ and define $\overline{t}_{k}$ as the ordered union of $t_k$ and $\tilde{t}_k$. 
    Then $\xi(\overline{t}_k)$, by a similar argument as before, converges. 
    This implies $\lim\limits_{k\to\infty}\xi(\tilde{t}_k)=\lim\limits_{k\to\infty}\xi(t_k)$. 
    Hence, the limit exists.
\end{beweis}

\begin{beweis}[ of \cref{thm:strong_convergence_Gradflow}]
    Short-time existence, long-time existence and subonvergence towards a critical point of the flow were shown in \cite[Theorem 1.7]{FrechesSchumacherSteenebruggevonderMoselPalaisSmaleConditionGeometric2025} and condition \textrm{L**} was verified in the proof of \cite[Theorem 1.3]{FrechesSchumacherSteenebruggevonderMoselPalaisSmaleConditionGeometric2025}. 
    By \cite[Corollary 3.13]{FrechesSchumacherSteenebruggevonderMoselPalaisSmaleConditionGeometric2025}, critical points of $ \TPp\vert_{\AL^s_{0}} $ are also critical for $ \TPp\vert_{\AL^s} $, thus 
    \[
        D^2\left( \TPp\vert_{\AL^s_{0}} \right)(\gamma_{\infty})=D^2\left( \TPp\vert_{\AL^s_{0}} \right)(\gamma_{\infty})\big\vert_{\T_{\gamma_{\infty}}\AL^s_0},
    \] 
    which is clearly is Fredholm with index zero as well. 
    Condition \textrm{L**} was shown in the proof of  \cite[Theorem 1.3]{FrechesSchumacherSteenebruggevonderMoselPalaisSmaleConditionGeometric2025}. 
    Thus, the claim follows from \cref{lem:finite_length_trajectories}, \cref{rem:convergence_gradflow_fixedpoint} and \cref{rem:LSineq_fixedpoint}. 
\end{beweis}

The \L{}ojasiewicz-Simon gradient inequality also provides information about the convergence rate of the gradient flow, see \cite[Theorem 24.21]{FeehanGlobalExistenceConvergence2016}. 
\begin{corollary}
	Let $\xi:[0,\infty)\to \AL^s$ and $\gamma_\infty$ as in \cref{convergence_gradflow}, then we have
	\begin{equation}
		\pdist_{\AL^s}(\xi(t),\gamma_\infty)\le C \Phi(g(t)),t\ge0
	\end{equation}
where
\begin{align*}
	\Phi(g(t))=\begin{cases}
		\frac{1}{Z(1-\theta)}\left(Z^2(2\theta-1)t+(\TPp(\gamma_0))^{1-2\theta}\right)^{-(1-\theta)/(2\theta-1)}, &\frac{1}{2}<\theta<1\\
		\frac{2}{Z}\sqrt{\TPp(\gamma_0)}\exp(-Z^2t/2),&\theta=\frac{1}{2}.
	\end{cases}
\end{align*}
\end{corollary}

As of now, it is unknown whether $\theta$ is optimal in our case or not. It is known however, that $\theta=\frac{1}{2}$ if and only if the energy is Morse-Bott, see \cite{FeehanMaridakisLojasiewiczSimonGradientInequalities2020}. We conjecture, that this is the case, and might return to this question in future projects.\\

\section*{Acknowledgments}
The authors want to thank Heiko von der Mosel and Philipp Reiter for their invaluable feedback during the preparation of this article.\\
Elias Döhrer gratefully acknowledges the funding support from the European Union and the Free State of Saxony (ESF).\\
Nicolas Freches is partially funded by {the DFG}-Graduiertenkolleg \emph{Energy, Entropy, and Dissipative Dynamics (EDDy)},  project no. 320021702/GRK2326.

\printbibliography

\end{document}